\documentclass[a4paper]{article}
\usepackage[utf8]{inputenc}

\usepackage[top=2.5cm, bottom=2.5cm, left=2.5cm, right=2.5cm]{geometry}
\usepackage{amsmath, amssymb, scrextend}
\usepackage[pdfencoding=auto, psdextra]{hyperref}

 \usepackage{xcolor}

 \usepackage{amsthm}
 \newtheorem{thm}{Theorem}[section]
 \newtheorem{lm}[thm]{Lemma}
 \newtheorem{res}[thm]{Result}

 \theoremstyle{definition}
 
 \newtheorem{rmk}[thm]{Remark}
 \newtheorem{df}[thm]{Definition}

 \newcommand{\FF}{\mathbb F}

 \newcommand{\vspan}[1]{\left \langle #1 \right \rangle}
 \newcommand{\set}[1]{ \left \{ #1 \right \} }
 \newcommand{\sett}[2]{ \left\{ #1 \, \, || \, \, #2 \right \} }

 \newcommand{\pg}{\textnormal{PG}}

 \newcommand{\rk}{\textnormal{rk}}
 \newcommand{\one}{\mathbf 1}
 \newcommand{\col}{\textnormal{col}}
 
 \newcommand{\mq}{\mathcal Q}

\title{Erd\H{o}s-Ko-Rado theorems for ovoidal circle geometries and polynomials over finite fields}
\author{Sam Adriaensen\thanks{Vrije Universiteit Brussel, Pleinlaan 2, 1050 Elsene, Belgium. Email: \url{sam.adriaensen@vub.be}  } \\ {\it Vrije Universiteit Brussel}}
\date{}

\begin{document}

\maketitle

\begin{abstract}
In this paper we investigate Erd\H{o}s-Ko-Rado theorems in ovoidal circle geometries.
We prove that in Möbius planes of even order greater than 2, and ovoidal Laguerre planes of odd order, the largest families of circles which pairwise intersect in at least one point, consist of all circles through a fixed point.
In ovoidal Laguerre planes of even order, a similar result holds, but there is one other type of largest family of pairwise intersecting circles.

As a corollary, we prove that the largest families of polynomials over $\FF_q$ of degree at most $k$, with $2 \leq k < q$, which pairwise take the same value on at least one point, consist of all polynomials $f$ of degree at most $k$ such that $f(x) = y$ for some fixed $x$ and $y$ in $\FF_q$.

We also discuss this problem for ovoidal Minkowski planes, and we investigate the largest families of circles pairwise intersecting in two points in circle geometries.
\end{abstract}

\bigskip


\section{Introduction}

In their seminal paper \cite{erdoskorado}, Erd\H{o}s, Ko, and Rado proved their famous theorem which states the following.

\begin{thm}[{\cite{erdoskorado}}]
Choose integers $k$ and $n$ such that $0 < k \leq n/2$.
Let $\mathcal F$ be a family of subsets of size $k$ of $\set{1,\dots,n}$ such that for each $F, G \in \mathcal F$, $F \cap G \neq \emptyset$.
Then
\[
 | \mathcal F | \leq \binom{n-1}{k-1}.
\]
Moreover, if $k < n/2$, then equality holds if and only if $\mathcal F$ is the family of all $k$-sets through a fixed element of $\set{1,\dots,n}$.
\end{thm}

This theorem inspired many generalisations.
Most of these can be formulated in terms of incidence structures.
Consider a set $\mathcal P$, whose elements are called points, and let $\mathcal B$ be a collection of $k$-element subsets of $\mathcal P$ for some constant $k$.
The elements of $\mathcal B$ are called blocks.
We call two blocks $B_1$ and $B_2$ \emph{$t$-intersecting} if $|B_1 \cap B_2| \geq t$.
We say that a family $\mathcal F \subseteq \mathcal B$ is $t$-intersecting if any two of its members are $t$-intersecting.
We call a 1-intersecting family simply an intersecting family.
Then we are interested in two questions:
\begin{enumerate}
 \item What size do the largest $t$-intersecting families in $(\mathcal P, \mathcal B)$ have?
 \item What structure do the largest $t$-intersecting families have?
\end{enumerate}

Given a $t$-element subset $S$ of $\mathcal P$, define $F(S) = \sett{B \in \mathcal B}{S \subseteq B}$.
Then $F(S)$ is clearly a $t$-intersecting family.
We say that the $t$-intersecting families in $(\mathcal P, \mathcal B)$ satisfy the \emph{weak EKR property} if the largest $t$-intersecting families have size $|F(S)|$ for some $t$-element set $S$, and that they satisfy the \emph{strong EKR property} if all $t$-intersecting families of maximum size are of the form $F(S)$.

\bigskip

The above theorem deals with the case where $\mathcal P = \set{1, \dots, n}$, $\mathcal B$ consists of all $k$-element subsets of $\mathcal P$, and $t=1$.
It states that intersecting families satisfy the weak EKR property if $k \leq n/2$, and the strong EKR property if $k < n/2$.
In the same paper, the authors prove that $t$-intersecting families in this incidence structure satisfy the strong EKR property if $n$ is large enough with respect to $k$ and $t$.
There have been a lot of improvements on finding the correct bound on $n$.
Eventually, Ahlswede and Kachatrian \cite{ahlswedekhachatrian} completed the study and characterised the largest $t$-intersecting families for all values of $n$, $k$, and $t$.

The $q$-analog of this problem can be formulated as finding the largest families $\mathcal F$ of $k$-dimensional subspaces of $\FF_q^n$, such that for all $F, G \in \mathcal F$, $\dim(F \cap G) \geq t$.
Hsieh \cite{hsieh75} proved that the weak EKR property holds for $t$-intersecting families if $n \geq 2k + 2$.
This was the first step to studying the Erd\H{o}s-Ko-Rado problem in finite vector spaces, or equivalently in finite projective spaces, which received a lot of attention since, see e.g.\ \cite{maartenthesis}.
Related structures such as the finite classical polar spaces have also been investigated, see e.g.\ \cite{pepestormevanhove}.

Another example of incidence structures where these questions have been posed are graphs of functions.
Consider two finite sets $A$ and $B$, and consider a subset $\Phi$ of functions $\varphi: A \to B$.
For each function $\varphi$, define its graph to be $\sett{(a,\varphi(a))}{a \in A}$.
Consider the incidence structure with $A \times B$ as points, and the graphs of the elements of $\Phi$ as blocks.
Then we can ask what the largest $t$-intersecting families are.
This problem can be reformulated without explicit mention of the graphs of the functions, by calling two functions $\varphi$ and $\psi$ $t$-intersecting when their graphs are $t$-intersecting, i.e.\ when there are at least $t$ values $a \in A$ for which $\varphi(a) = \psi(a)$.
Frankl and Deza \cite{frankldeza} proved that intersecting families in the symmetric group satisfy the weak EKR property.
Cameron and Ku \cite{cameronku} proved that they satisfy the strong EKR property.
One of the most general results in this context is due to Meagher, Spiga and Tiep \cite{meagherspigatiep}, who proved that intersecting families in any faithful 2-transitive permutation group satisfy the weak EKR property.

These are just a few of the numerous results inspired by Erd\H{o}s, Ko, and Rado.
We refer to interested reader to \cite{godsilmeagher}, which to author cannot but recommend as an excellent introduction into this subject and into the algebraic techniques that will be used in this paper.

\bigskip

Now we can formulate the main theorems of this paper.
The first one concerns polynomials of bounded degree over finite fields.

\begin{thm}
 \label{ThmMainPolynomials}
Assume that $t \leq k < q$.
Consider a $t$-intersecting family $\mathcal F$ of polynomials over $\FF_q$ such that the degree of each polynomial in $\mathcal F$ is at most $k$.
Then $|\mathcal F| \leq q^{k+1-t}$.
Furthermore, if $t=1$ and $k\geq 2$, then equality holds if and only if there are elements $x,y \in \FF_q$ such that
\[
 \mathcal F = \sett{f \in \FF_q[X]}{\deg f \leq k, \; f(x)=y}.
\]
\end{thm}

In other words, $t$-intersecting families of polynomials of degree at most $k$ satisfy the weak EKR property.
If $k \geq 2$ and $t=1$, than they satisfy the strong EKR property.

\bigskip

The other main results of this paper concern the so-called circle geometries.
There are three types of circle geometries, the Möbius planes, Laguerre planes, and Minkowski planes.
We refer to \S \ref{SubsectionCircleGeometries} for the definitions.
These are certain types of incidence structures in which the blocks are often referred to as circles.

\begin{thm}
 \label{ThmMainMobius}
Let $M = (\mathcal P, \mathcal B)$ be a Möbius plane of even order $q>2$.
Let $\mathcal F \subseteq \mathcal B$ be an intersecting family.
Then $|\mathcal F| \leq q(q+1)$.
Equality occurs if and only if $\mathcal F$ consists of all circles through a point of $\mathcal P$.
\end{thm}

\begin{thm}
 \label{ThmMainLaguerreOdd}
Let $L = (\mathcal P, \mathcal B)$ be an ovoidal Laguerre plane of order $q$, $q$ odd.
Let $\mathcal F \subseteq \mathcal B$ be an intersecting family.
Then $|\mathcal F| \leq q^2$.
Equality occurs if and only if $\mathcal F$ consists of all circles through a point of $\mathcal P$.
\end{thm}

\begin{thm}
 \label{ThmMainLaguerreEven}
Let $L = (\mathcal P, \mathcal B)$ be an ovoidal Laguerre plane of order $q$, $q>2$ even.
Let $\mathcal F \subseteq \mathcal B$ be an intersecting family.
Then $|\mathcal F| \leq q^2$.
Equality occurs if and only if $\mathcal F$ consists of all circles through a point of $\mathcal P$ or $\mathcal F$ consists of a circle $c$ in $\mathcal B$ and all circles intersecting $c$ in exactly one point.
\end{thm}

In other words, intersecting families in Möbius planes of even order greater than 2 and ovoidal Laguerre planes of odd order satisfy the strong EKR property.
The same holds in ovoidal Minkowski planes, which follows from a result by Meagher and Spiga \cite{meagherspiga}, see Theorem \ref{ThmMinkowski}.
On the other hand, intersecting families in ovoidal Laguerre planes of even order satisfy the weak EKR-property, but not the strong EKR-property.

\bigskip

We also give upper bounds on the size of 2-intersecting families in Laguerre and Minkowski planes.
We use the same arguments as Blokhuis and Bruen \cite{blokhuisbruen} used for the Möbius planes.
The general bounds are approximately $\frac{q^2} 2$.
However, computations show that for small values of $q$, the size of a 2-intersecting family doesn't exceed $2q$.
We only find a satisfactory upper bound on the size of 2-intersecting families in ovoidal Laguerre planes of even order, where the weak EKR property is satisfied.

\bigskip

In \S 2, we give the necessary preliminaries, in \S 3 we give an overview of the rest of the paper.

\section{Preliminaries}

\subsection{Notation}

Throughout this article, $q$ denotes a prime power, $\FF_q$ denotes the finite field of order $q$, and $\FF_q^*$ denotes $\FF_q \setminus \set 0$.
If $V_1, \dots, V_n$ is a sequence of subsets of some vector space or projective space, their span is denoted by $\vspan{V_1,\dots,V_n}$.
The $n \times n$ all-one matrix is denoted as $J_n$, the all-one vector of length $n$ is denoted as $\one_n$, the $n \times n$ identity matrix is denoted as $I_n$, and we drop the indices if they are clear from context.
If the rows and columns of a matrix $A$ are indexed by sets $X$ and $Y$ respectively, then we will denote the entry in $A$ in the row corresponding to $x \in X$ and in the column corresponding to $y \in Y$ by $A(x,y)$.
The smallest eigenvalue of $A$ is denoted by $\tau(A)$.
The corresponding eigenspace is called the $\tau$-eigenspace.
The column space of $A$ is denoted by $\col(A)$.
The projective space of (projective) dimension $n$ over $\FF_q$ is denoted as $\pg(n,q)$.
Given a point set $\mq$ in $\pg(n,q)$, we call a line \emph{skew}, \emph{tangent}, or \emph{secant} to $\mq$ if it intersects $\mq$ in 0, 1, or 2 points respectively.
More generally, we call a line an $i$-secant if it intersects $\mq$ in exactly $i$ points.

\subsection{Quadratic sets in projective spaces}
 \label{SubsectionQuadraticSets}

A \emph{quadric} in $\pg(n,q)$ is the set of points $(x_0, \dots, x_n)$ in $\pg(n,q)$ which satisfy some homogeneous equation of degree two $\sum_{i,j=0}^n a_{i,j} X_i X_j = 0$, for some $a_{i,j} \in \FF_q$.
A quadric is called \emph{non-degenerate} if its equation is irreducible, or equivalently if it has no \emph{singular} point.
A point is singular if it only lies on tangent lines and lines contained in the quadric.
A non-degenerate quadric in $\pg(2,q)$ is called a \emph{conic}.

\bigskip

An \emph{oval} $\mathcal O$ in $\pg(2,q)$ is a set of $q+1$ points, no three collinear.
One of the most fundamental results in finite geometry states that if $q$ is odd, the only ovals in $\pg(2,q)$ are the conics.
This result was proven by Segre \cite{segre55}.

If $q > 8$ is even, then other examples of ovals in $\pg(2,q)$ are known.
However, in this case, every oval $\mathcal O$ has a \emph{nucleus} \cite{segre55}.
This is a point $N \notin \mathcal O$ such that the tangent lines to $\mathcal O$ are exactly the lines through $N$.
It readily follows that every point $P \notin \mathcal O \cup \set N$ lies on $\frac q 2$ skew lines, one tangent line, and $\frac q 2$ secant lines to $\mathcal O$.

For both $q$ odd and even, there are $\binom {q+1} 2$ secant lines, $q+1$ tangent lines, and $\binom q 2$ skew lines to $\mathcal O$.
If $q$ is odd, the tangent lines to $\mathcal O$ form a dual conic.
Therefore, the points not in $\mathcal O$ consist of $\binom{q+1}2$ so-called \emph{external} points, which lie on 2 tangent and $\frac{q-1}2$ secant and skew lines, and $\binom q 2$ \emph{internal} points, which lie on no tangents, and $\frac{q+1}2$ secant and skew lines.
Dually, a secant line contains $\frac{q-1}2$ external and internal point, a skew line contains $\frac{q+1}2$ external and internal points.

We also remark that up to isomorphism, there is only one conic in $\pg(2,q)$, which has equation $X_0 X_2 = X_1^2$.

\bigskip

Buekenhout \cite{buekenhout69} defined the more general concept of a \emph{quadratic set} as a set $\mathcal Q$ in $\pg(n,q)$ satisfying the following properties.
\begin{enumerate}
 \item Any line that intersects $\mathcal Q$ in more than two points is contained in $\mathcal Q$.
 \item Given a point $P \in \mathcal Q$, the union $\mathcal Q_P$ of the lines $l$ through $P$ for which $l \cap \mathcal Q$ is either $\set P$ or $l$, is either a hyperplane (called the \emph{tangent} hyperplane to $P$) or the entire space.
\end{enumerate}
If $\mathcal Q_P$ is the entire space for some point $P \in \mathcal Q$, we call $P$ a \emph{singular} point.
If there are no singular points, then $\mathcal Q$ is a \emph{non-degenerate} quadratic set.

Now suppose that $\mathcal Q$ is a quadratic set in $\pg(3,q)$, which is not the union of two (not necessarily distinct) subspaces.
Then there are three options:

\begin{enumerate}

\item First assume that $\mathcal Q$ is non-degenerate and does not contain a line.
Barlotti \cite[\S 5.2]{barlotti65} proved that if $q>2$, than a set of points, no three of which are collinear, contains at most $q^2+1$ points.
Sets of size $q^2+1$ with no triple of collinear points are called \emph{ovoids}.
Barlotti also showed that if $q>2$, every plane intersects an ovoid in either a point or an oval, proving that all ovoids are quadratic sets.
Extending Segre's result, Barlotti \cite{barlotti55} proved that if $q$ is odd, every ovoid of $\pg(3,q)$ is a so-called elliptic quadric.
Up to isomorphism $\pg(3,q)$ contains only one elliptic quadric, which has equation $g(X_0,X_1) + X_2 X_3 = 0$ for some irreducible quadratic form $g$.

The only known ovoids which are not elliptic quadrics are the Suzuki-Tits ovoids, constructed in \cite{tits60}, which exist for $q=2^h$, $h>1$ odd.

If $q=2$, then the complement of a plane is a collection of $q^3$ points, no three collinear, which is the reason why we make the distinction.

\item Now assume that $\mathcal Q$ is non-degenerate and $\mathcal Q$ does contain a line.
Buekenhout \cite{buekenhout69} proved that any non-degenerate quadratic set in $\pg(n,q)$, $n \geq 3$, which contains a line, is a non-degenerate quadric.
In $\pg(3,q)$ the only such quadrics are the hyperbolic quadrics, which are isomorphic to the quadric with equation $X_0 X_1 + X_2 X_3 = 0$.

\item Lastly, assume that $\mathcal Q$ is degenerate.
Let $R$ be a singular point.
Thus, the lines through $R$ which contain another point of $\mathcal Q$ are completely contained in $\mathcal Q$.
Take a plane $\pi$ not through $R$.
Let $\mathcal O$ be the intersection of $\pi$ and $\mathcal Q$.
Then it is easy to check that $\mathcal Q$ is the \emph{cone} with \emph{base} $\mathcal O$ and \emph{vertex} $R$, i.e.\ $\mathcal Q = \bigcup_{P \in \mathcal O} \vspan{R,P}$.
The fact that $\mathcal Q$ is a quadratic set implies that $\mathcal O$ is a quadratic set.
The only option is that $\mathcal O$ is an oval, otherwise $\mathcal Q$ is the union of two subspaces.
We will refer to this type of quadratic set as an \emph{oval cone}.
Moreover, if the oval $\mathcal O$ is a conic, then $\mathcal Q$ is a quadric, to which we will refer as a \emph{quadric cone}.
\end{enumerate}

\subsection{Circle geometries}
 \label{SubsectionCircleGeometries}

There are three types of circle geometries, also referred to as Benz planes, due to Benz's large contribution to the theory \cite{benz73}.
We give a description of these circle geometries, which can be found in \cite{hartmann}.
In this paper, we are only interested in the finite case.
The three geometries can be defined in a unified fashion.
Every circle geometry is an incidence structure $(\mathcal P, \mathcal B)$.
The elements of $\mathcal B$ are called circles.
A set of points is called \emph{concyclical} if they are contained in a circle.
A \emph{parallel relation} is a partitioning of the points $\mathcal P$ into \emph{parallel classes}.
We call two points parallel (with respect to a certain parallel relation) if they belong to the same parallel class.
Now assume that any two distinct points are parallel w.r.t.\ at most one parallel relation.
For every point $P \in \mathcal P$, we define the \emph{residue} at $P$ as the incidence structure $(\mathcal P_P, \mathcal B_P)$, where $\mathcal P_P$ are the points which are not parallel to $P$.
To define $\mathcal B_P$, let $\mathcal B_P'$ be the set containing all elements of $\mathcal B$ through $P$ and all parallel classes not containing $P$.
Then $\mathcal B_P = \sett{B \cap \mathcal P_P}{B \in \mathcal B_P'}$.

\begin{df}
 A \emph{circle geometry} is an incidence structure with at most two parallel relations, such that the residue at each point is an affine plane.
 We call a circle geometry a \emph{Möbius plane}, \emph{Laguerre plane}, or a \emph{Minkowski plane} if it has respectively 0, 1, or 2 parallel relations.
\end{df}

The characteristic properties of these circle geometries are that
\begin{enumerate}
 \item every three pairwise non-parallel points lie on a unique circle,
 \item given two non-parallel points $P$ and $Q$ and a circle $c$ through $P$ but not through $Q$, there is a unique circle through $P$ and $Q$ which intersects $c$ only in $P$,
 \item given a point $P$ and a circle $c$ not through $P$, $c$ intersects every parallel class through $P$ in a unique point.
\end{enumerate}

For each point $P$ of a circle geometry, the order of the affine plane $(\mathcal P_P, \mathcal B_P)$ is equal.
If this number is $q$, we say that the circle geometry has order $q$.
This means that every circle contains $q+1$ points, every point lies on $q(q+1-r)$ circles, where $r$ denotes the number of parallel relations, and every parallel class contains $q+r-1$ points.

A Möbius plane of order $q$ has $q^2+1$ points and $q(q^2+1)$ circles.
A Laguerre plane of order $q$ has $q(q+1)$ points and $q^3$ circles.
A Minkowski plane of order $q$ has $(q+1)^2$ points and $(q-1)q(q+1)$ circles.

\bigskip

The quadratic sets described above can be used to construct the \emph{ovoidal} circle geometries.
Let $\mathcal Q$ be a quadratic set in $\pg(3,q)$, which is not the union of two subspaces.
Call a plane $\pi$ an \emph{oval plane} if $\mq$ intersects $\pi$ in an oval.
Consider the incidence structure $(\mathcal P, \mathcal B)$, where $\mathcal P$ are the non-singular points of $\mathcal Q$, and $\mathcal B$ are the intersections of $\mathcal Q$ with the oval planes.
If $\mathcal Q$ is an ovoid, oval cone, or hyperbolic quadric, then the circle geometry is respectively a Möbius plane, Laguerre plane, or Minkowski plane.

There is an extra axiom for circle geometries, called the Bundle Theorem, which states the following.

\bigskip

{\bf Bundle Theorem.}
Consider 8 pairwise non-parallel points $A_1,A_2,A_3,A_4,B_1,B_2,B_3,B_4$.
Define for $i<j$ the set $Q_{i,j} = \set{A_i,A_j,B_i,B_j}$.
Assume that for 5 of the 6 choices for $(i,j)$, the four points of $Q_{i,j}$ are concyclical, and at least 4 of these circles are distinct.
Then for the sixth choice of $(i,j)$, the points of $Q_{i,j}$ are concyclical as well.

\bigskip

Kahn \cite{kahn80} proved that a circle geometry is ovoidal if and only if it satisfies the Bundle Theorem.
Often this result is compared the fact that a projective plane satisfies the Theorem of Desargues if and only if it arises from a vector space over a skewfield.

Moreover, there is a stronger axiom, which plays a similar role.

\bigskip

{\bf Theorem of Miquel.}
Take eight points $P_1, \dots, P_8$ and assign them to the vertices of a cube.
Assume that for 5 faces of the cube, its 4 vertices are concyclical.
Then the same holds for the 6\textsuperscript{th} face of the cube.

\bigskip 

Chen \cite{chen70,chen74} proved that a circle geometry satisfies the Theorem of Miquel if and only if it is ovoidal, and the quadratic set used to construct the geometry is a quadric.
This can be seen as an analog of the Theorem of Pappus for projective planes, which holds in a projective plane if and only if it arises from a vector space over a (commutative) field.
Furthermore, by the results of Segre, Barlotti, and Buekenhout, we know that in finite circle geometries of odd order, and Minkowski planes of even order, the Bundle Theorem is equivalent with the Theorem of Miquel.
This result can be seen as an analog of Wedderburn's little Theorem, which implies that the Theorems of Desargues and Pappus are equivalent for finite projective planes.
The pièce de résistance is that there is an analog of the infamous prime power conjecture for finite projective planes, which remains perhaps the most elusive open problem in finite geometry.

\begin{res}[{\cite{dembowski64}}]
 \label{DoYouEvenMobius}
Any finite Möbius plane of even order is ovoidal.
\end{res}

The interested reader is also refered to \cite[\S 6]{dembowksi68} for a treatise on Möbius planes.
We remark that Möbius planes appear there under the name ``inversive planes''.

\subsection{Association schemes}

Many tools have been developed to tackle EKR problems, the strongest of which might be the theory of association schemes, which we revise here.
The most important contributions were made by Delsarte \cite{delsartethesis}.
The interested reader is also referred to the surveys \cite{godsil}, \cite[\S 3]{godsilmeagher}, \cite[\S 2]{bcn}.

\begin{df}
 \label{DefAssoc1}
Given a set $X$, a \emph{(symmetric) d-class association scheme} on $X$ is a partition $R_0, \dots, R_d$ of $X \times X$ such that 
\begin{enumerate}
 \item $R_0 = \sett{(x,x)}{x \in X}$,
 \item for each $x,y \in X$ and each $i = 0, \dots, d$, if $(x,y) \in R_i$, then $(y,x) \in R_i$,
 \item given $(x,y) \in R_k$, the number of elements $z \in X$ such that $(x,z) \in R_i$ and $(z,y) \in R_j$ is some number $p^k_{i,j}$, which is independent of $(x,y)$.
\end{enumerate}
\end{df}

The numbers $p_{i,j}^k$ are called the \emph{intersection numbers} of the association scheme.
Note that condition (2) and (3) imply that $p_{i,j}^k = p_{j,i}^k$.

Ofttimes, the definition is reformulated in terms of \emph{adjacency matrices}.
With each relation $R_i$, one can associate a matrix $A_i$ such that the rows and columns of $A$ are indexed by $X$, and
\[
 A_i(x,y) = \begin{cases}
 1 & \text{if } (x,y) \in R_i, \\
 0 & \text{otherwise.}
 \end{cases}
\] 

Then the conditions imposed on $R_0, \dots, R_d$ in order to be an association scheme, translate to the following equivalent conditions on the matrices $A_0, \dots, A_d$.
\begin{enumerate}
 \item $A_0 = I$,
 \item for each $i=0,\dots,d$, $A_i^t = A_i$,
 \item $\displaystyle A_i A_j = \sum_{k=0}^d p_{i,j}^k A_k$.
\end{enumerate}

Note that conditions (2) and (3) imply that $A_i$ and $A_j$ commute.

The matrices $A_0, \dots, A_d$ are also called an association scheme, since they convey the same information as $R_0, \dots, R_d$.
The fact that $R_0, \dots, R_d$ is a partition of $X \times X$ is equivalent with $A_0 + \ldots + A_d = J$.

Since $A_0, \dots, A_d$ are symmetric and commute, they can be diagonalised simultaneously.
The association scheme encodes important information about this simultaneous diagonalisation.
We call a matrix $E$ an \emph{orthogonal projection} if $E^2 = E$, and $\col(E) = \ker(E)^\perp$.

\begin{thm}
Given a $d$-class association scheme $A_0, \dots, A_d$ of $n \times n$-matrices, there exist a basis $E_0 = \frac 1 n J_n$, $E_1, \dots, E_d$ of orthogonal projections of $\vspan{A_0, \dots, A_d}$ such that $E_0 + \dots + E_d = I_n$.
\end{thm}

A direct corollary of this theorem is that if we define $V_i = \col(E_i)$, then $V_i$ are pairwise orthogonal subspaces, spanning the entire space.
Every eigenspace of each matrix in the association scheme is a direct sum of some of these spaces.
The spaces $V_i$ will be referred to as the \emph{eigenspaces of the association scheme}.

Now suppose that $A_0, \dots, A_d$ is a $d$-class association scheme on a set $X$ with eigenspaces $V_0, \dots, V_d$.
Let $P$ be the $(d+1)\times(d+1)$-matrix indexed by $\set{0,\dots,d}$, such that $P_{i,j}$ is the eigenvalue of $A_j$ on the eigenspace $V_i$.
$P$ encapsulates the essential information of the association scheme.
We call $P$ the \emph{matrix of eigenvalues} of the association scheme, and $Q := |X| P^{-1}$ the \emph{dual matrix of eigenvalues}.
It is well-known that the first row of $Q$ is $(\dim V_0, \dots, \dim V_d)$.

\bigskip 

Throughout this article, given an association scheme $R_0, \dots, R_d$, we will use the notation $A_i$ and $V_i$ for the adjacency matrices and eigenspaces as introduced above.

\subsection{Bounds on cocliques}
 \label{SUbsectionBoundsCocliques}

Suppose that $A$ is a symmetric $01$-matrix whose rows and columns are indexed by the same set $X$, and whose diagonal is zero.
Then we can associate to $A$ the graph $G$ whose vertices are the elements of $X$ and where two vertices $x$ and $y$ are adjacent if and only if $A(x,y) = 1$.
We call $A$ the \emph{adjacency matrix} of $G$.
A \emph{clique} in $G$ is a set of pairwise adjacent vertices, a \emph{coclique} is a set of pairwise non-adjacent vertices.
Let $\omega(G)$ and $\alpha(G)$ denote the size of the largest clique and coclique respectively of the graph $G$.
The \emph{characteristic vector} of a subset $K$ of the vertices of $G$ is the $01$-vector indexed by the vertices of $G$, which has 1 exactly in the positions corresponding to the elements of $K$.
The following bound is known as the ratio bound or as Hoffman's bound.
As explained in a recent note by Haemers \cite{haemers2021}, it is tricky to give a correct reference for this bound.
A proof can be found e.g.\ in \cite[Theorem 2.4.2]{godsilmeagher}.

\begin{thm}[Hoffman's ratio Bound]
 \label{ThmRatioBound}
Let $G$ be a $k$-regular graph on $n$ vertices with adjacency matrix $A$.
Denote $\tau = \tau(A)$.
Then
\[
 \alpha(G) \leq \frac n {1 + \frac k {-\tau}}.
\]
Moreover, if $K$ is a coclique attaining this bound, then the characteristic vector of $K$ lies in the sum of the $k$- and $\tau$-eigenspaces of $A$.
\end{thm}

There is another interesting bound on the size of cliques and cocliques.
This bound is due to Delsarte \cite[Theorem 3.9]{delsartethesis} in the context of association schemes, but holds for some other classes of graphs as well (see e.g.\ \cite[Corollary 2.1.2]{godsilmeagher} for a proof).

\begin{thm}[Clique-coclique bound]
If $G$ is a vertex-transitive graph on $n$ vertices, then $\alpha(G) \omega(G) \leq n$.
\end{thm}

\subsection{Intersection matrices of an association scheme}
\label{SubsectionIntersectionMatrices}

Let $R_0, \dots, R_d$ be the relations of an association scheme on a set $X$ of size $n$.
Let $p^k_{i,j}$ be the intersection numbers.
Let $A_i$ denote the adjacency matrix of $R_i$, and $G_i$ the corresponding graph.
Now take a vertex $x \in X$.
Partition $X$ in the classes $T_i = \sett{y \in X}{(x,y) \in R_i}$.
Let $T_x$ denote the $n \times (d+1)$-matrix whose columns are the characteristic vectors of $T_0, \dots, T_d$.
Consider the matrix $A_i T$.
The rows of $A_i T$ correspond to elements of $X$ and the columns to the classes $T_j$.
The entry in position $(y,T_j)$ equals the number of elements $z \in X$ such that $(y,z) \in R_i$ and $(x,z) \in R_j$.
Therefore, this number equals $p^k_{i,j}$ if $(x,y) \in R_k$, or equivalently if $y \in T_k$.
Hence, if $B_i$ denotes the matrix $(p^k_{i,j})_{k,j}$, then $A_i T_x = T_x B_i$.
In the words of graph theory, $T_0, \dots, T_d$ are an \emph{equitable partition} of each graph $G_i$.

As described in \cite[Proposition 2.2.2]{bcn}, the matrices $B_0,\dots,B_d$, called the \emph{intersection matrices}, can be simultaneously diagonalised by the columns of the matrix $Q$ of dual eigenvalues.
Furthermore, if $v$ is an eigenvector for $B_i$ with eigenvalue $\lambda$, then $A_i T_x v = T_x B_i v = \lambda T_x v$.
Hence, $T_x v$ is an eigenvector of $A_i$ with eigenvalue $\lambda$.
We call $T_x v$ a \emph{lift} of the eigenvector $v$.
For each column $Q_j$ of $Q$ and for each $x \in X$, $T_x Q_j$ is an eigenvector which lies in the eigenspace $V_j$.
Therefore, $Q^{-1} B_i Q$ is a diagonal matrix with as diagonal the $i$\textsuperscript{th} column of $P$.

Since $Q$ is invertible, $Q (w_0, \dots, w_d)^t = (1,0,\dots,0)^t$ for some vector $(w_0, \dots, w_d)$.
Therefore, $$w_0 T_x Q_0 + \dots + w_d T_x Q_d = T_x Q \begin{pmatrix} w_0 \\ \vdots \\ w_d \end{pmatrix} =  T_x \begin{pmatrix} 1 \\ 0  \\ \vdots \\ 0 \end{pmatrix}.$$
Note that $T_x (1,0,\dots,0)^t$ is the unit vector with a 1 in the position corresponding to $x$.
Thus, all the lifts of all the columns of $Q$ generate $\mathbb C^n$.
As a consequence, all the lifts of $Q_j$ must generate $V_j$.

We can use the intersection matrices in this way to determine the matrix $P$, and lift the columns of $Q$ to find generating sets of the eigenspaces of the association scheme.

\subsection{Kronecker products}
 \label{SubsectionKronecker}

Given two matrices $A$ and $B$, we denote their Kronecker product by $A \otimes B$.
This means that $A \otimes B$ is the matrix obtained by replacing every entry $a_{i,j}$ of $A$ by the matrix $a_{i,j} B$.
We use the following property of the Kronecker product, see e.g.\ \cite[page 65]{willihans}.

\begin{lm}
 \label{LmKronecker}
For any matrices $A,B,C,D$ for which the following products are well-defined,
$(A \otimes B) (C \otimes D) = (A C) \otimes (B D)$.
\end{lm}

Suppose that $A$ and $B$ are square matrices, that $v_1 \dots, v_n$ is a basis of eigenvectors of $A$ with $A v_i = \lambda_i v_i$, and $w_1, \dots, w_m$ is a basis of eigenvectors for $B$ with $B w_j = \mu_j \lambda_j$.
Then $(A \otimes B) (v_i \otimes w_j) = (A v_i) \otimes (B w_j) = (\lambda_i \mu_j) (v_i \otimes w_j)$.
Furthermore, if $V$ is the matrix with columns $v_1, \dots, v_n$ and $W$ is the matrix with columns $w_1, \dots, w_m$,
then $(V^{-1} \otimes W^{-1}) (v_i \otimes w_j) = e_i \otimes e_j$, where $e_i$ and $e_j$ denote the standard unit vectors.
Since the vectors $e_i \otimes e_j$ are a basis, the vectors $v_i \otimes w_j$ are a basis of eigenvectors of $A \otimes B$.

\section{Overview}

As mentioned in the introduction, the main goal of this paper is proving that a large class of circle geometries satisfy the strong EKR property for intersecting families.
This will be done in three different settings, but the proofs will follow the same outline.
First we define the following relations.

\begin{df}
 \label{DfRelAssoc}
Let $(\mathcal P, \mathcal B)$ be a circle geometry.
We define the following relations on $\mathcal B$.
\begin{itemize}
 \item $R_0 = \sett{(c,c)}{c \in \mathcal B}$,
 \item $R_1 = \sett{(c,c') \in \mathcal B^2}{|c \cap c'| = 1}$,
 \item $R_2 = \sett{(c,c') \in \mathcal B^2}{|c \cap c'| = 2}$,
 \item $R_3 = \sett{(c,c') \in \mathcal B^2}{|c \cap c'| = 0}$.
\end{itemize}
\end{df}
We will prove that these constitute an association scheme for some classes of circle geometries.
This will enable us to use the machinery of algebraic graph theory.
The relations are evidently symmetric and a partition of $\mathcal B \times \mathcal B$, so we will only need to prove that the intersection numbers are well-defined.
Note that if $p^k_{i,j}$ is well-defined, by the symmetry $p^k_{i,j} = p^k_{j,i}$.

Let $A_i$ denote the adjacency matrix of $R_i$.
Then a coclique in the graph $G$ with adjacency matrix $A_3$ is exactly the same thing as an intersecting family in the circle geometry.
We will use Hoffman's ratio bound to determine the maximum size of a coclique in $G$, and to determine the subspace in which the characteristic vectors of the maximum cocliques lie.
We will prove that the characteristic vectors of the largest intersecting families lie in the column space of the so-called incidence matrix.

\begin{df}
 \label{DfIncidenceMatrix}
The \emph{incidence matrix} of an incidence structure $(\mathcal P, \mathcal B)$ is the 01-matrix, whose rows are indexed by the elements of $\mathcal B$, whose columns are indexed by the elements of $\mathcal P$, and whose entry in position $(B,P)$ is 1 if $P \in B$, and 0 otherwise.
\end{df}

Then we need to prove that the only characteristic vectors of intersecting families of the appropriate size in the column space of the incidence matrix, are the the columns of the matrix itself.

\bigskip

In \S \ref{SectionMobius}, we use this strategy to prove Theorem \ref{ThmMainMobius} which states that intersecting families in Möbius planes of even order greater than two satisfy the strong EKR property.
In \S \ref{SectionLaguerre} we use this strategy to prove Theorem \ref{ThmMainLaguerreOdd} and \ref{ThmMainLaguerreEven}, and classify all the largest intersecting families in ovoidal Laguerre planes.
We also give bounds on 2-intersecting families in Laguerre planes.
For $q$ even, we find a sharp bound for the maximal size of a 2-intersecting family.
In \S \ref{SectionPolynomials}, we use our results on ovoidal Laguerre planes to prove that intersecting families in the set of polynomials over $\FF_q$ of degree at most $k \in [2,q-1]$ satisfy the strong EKR-property, i.e.\ Theorem \ref{ThmMainPolynomials}.
In \S \ref{SectionMinkowski} we briefly discuss these problems for Minkowski planes.
We make some concluding remarks in \S \ref{SectionGeMoogtNaarHuisGaan}.

\section{Möbius planes of even order}
 \label{SectionMobius}

Let $M = (\mathcal P, \mathcal B)$ be a Möbius plane of even order $q > 2$.
By Result \ref{DoYouEvenMobius}, $M$ is ovoidal.
Let $\mq$ be the ovoid in $\pg(3,q)$ from which $M$ is constructed.
Then $\mathcal P = \mq$ as point set, and $\mathcal B$ corresponds to the oval planes.
There are $q^3 + q$ oval planes, since there are $q^3 + q^2 + q + 1$ planes in $\pg(3,q)$, of which $|\mathcal P| = q^2 + 1$ are tangent planes.
If $\pi$ is an oval plane, the $\pi \cap \mq$ has a nucleus $N$.
We say that $N$ is the nucleus of $\pi$.

\begin{lm}
 Every point $P \notin \mq$ is the nucleus of a unique oval plane.
\end{lm}

\begin{proof}
Take a point $P \notin \mq$, and let $x$ denote the number of tangent planes through $P$.
Every point of $\mq$ lies on $q+1$ planes through $P$, and the total number of planes through $P$ equals $q^2+q+1$, whence
\[
 x + (q^2+q+1-x)(q+1) = (q^2+1)(q+1).
\]
This implies that $P$ lies in $x=q+1$ tangent planes.
Since each such tangent plane contains a unique point of $\mq$, $P$ lies on $q+1$ tangent lines.
The plane $\pi$ spanned by two such tangent lines must be an oval plane, and $P$ lies on two tangent lines to this oval.
Therefore, $P$ is the nucleus of this oval, and the tangent lines through $P$ are exactly the lines in $\pi$ through $P$.
\end{proof}

An easy counting argument shows that a skew, tangent, and secant line lies on respectively 2, 1, and 0 tangent planes.

\subsection{The association scheme}

We will now show that the relations from Definition \ref{DfRelAssoc} give us an association scheme.
We will do this by proving that the intersection numbers are well-defined.
First consider the numbers $p^0_{i,j}$.
Given an oval plane $\pi$, this is the number of oval planes $\rho$ such that $(\pi,\rho) \in R_i$ and $(\pi,\rho)\in R_j$.
Obviously, this is zero if $i \neq j$.
\begin{itemize}
 \item If $i=j=0$, then the only such plane $\rho$ is $\pi$ itself, hence $p^0_{0,0} = 1$.
 \item Now suppose that $i=j=1$.
There are $q+1$ tangent lines in $\pi$, which all lie on a unique tangent plane, thus on $q-1$ oval planes $\rho \neq \pi$.
Therefore, $p^0_{1,1} = (q+1)(q-1)$.
 \item Assume that $i=j=2$.
 There are $\binom{q+1}2$ secants to $\mq$ in $\pi$, each lying on no tangent planes, thus $p^0_{2,2} = (q+1) \frac{q^2}2$.
 \item Let $i$ and $j$ be $3$.
 There are $\binom{q}2$ skew lines to $\mq$ in $\pi$, each lying on 2 tangent planes, hence $p^0_{3,3} = q(q-1) \frac{q-2}2$.
\end{itemize}

Take two oval planes $\pi_1$ and $\pi_2$.
Let $N_1$ and $N_2$ denote their respective nuclei.
Denote $l = \pi_1 \cap \pi_2$.
Then for each $i$ and $j$ in $\set{0,1,2}$, we need to determine the number of oval planes $\pi$ that intersect $\pi_1$ in and $i$-secant and $\pi_2$ in a $j$-secant.
Such a plane $\pi$ either intersects $\pi_1$ and $\pi_2$ in $l$, or in distinct lines, which intersect in a point of $l$.
Vice versa, given a point $P$ of $l$, and lines $l_1$ in $\pi_1$ and $l_2$ in $\pi_2$ through $P$, the lines $l_1$ and $l_2$ span a unique plane.
If this plane is an oval plane, its relation to $\pi_1$ and $\pi_2$ only depends on whether $l_1$ and $l_2$ are an skew, tangent or secant line to $\mq$.
Therefore, we can compute the intersection numbers by counting the number of $i$- and $j$-secants in $\pi_1$ and $\pi_2$ which intersect, and taking into account that we need to disregard the tangent planes.

First suppose that $l \cap \mq$ is a single point $Q$.
Then $N_1$ and $N_2$ are distinct points of $l$.

We compute the numbers $p_{i,j}^1$.
\begin{itemize}
 \item $p^1_{1,1} = q-2$: an oval plane intersects $\pi_1$ and $\pi_2$ in one point of $\mq$ if and only if it contains $N_1$ and $N_2$, hence if and only if it contains $l$.
 The line $l$ lies on $q$ oval planes, two of which are $\pi_1$ and $\pi_2$.
 \item $p_{1,2}^1 = \frac{q^2}2$: A plane intersects $\pi_1$ in one point and $\pi_2$ in two points if and only if it intersects $\pi_1$ in a tangent line (through $N_1$), distinct from $l$, which gives $q$ options, and $\pi_2$ in secant line through $N_1$ which gives $\frac q 2$ options.
 \item $p_{1,3}^1 = q \frac{q-2}2$: There are $q$ ways to choose a tangent line $t \neq l$ in $\pi_1$ (through $N_1$).
 Then $t \cap \mq = P$.
 Through $N_1$ there are $\frac q 2$ skew lines in $\pi_2$.
 One of these skew lines is excluded from our count because it is the intersection of the tangent plane through $P$ and $\pi_2$.
 \item $p^1_{2,2} = \frac{q^2}4(q+2)$: We count the number of secant lines in $\pi_1$ and $\pi_2$ that intersect $l$ in the same point.
 Through each of the $q-2$ points of $l \setminus \set{N_1,N_2,Q}$ there are $\frac q 2$ secant lines in $\pi_1$ and $\pi_2$, through $Q$ there are $q$ secant lines in either plane, through $N_1$ and $N_2$ there are none.
 Thus, $p^1_{2,2} = (q-2)(q/2)^2 + q^2$.
 \item $p_{2,3}^1 = \frac{q^2}4(q-2)$: Analogous to $p_{2,2}^1$, except there are no skew lines through $Q$.
 \item $p_{3,3}^1 = \frac q 4 (q-2)(q-4)$: There are $q-2$ points on $l \setminus \set{Q,N_1,N_2}$, each on $\frac q 2$ skew lines in both planes, thus $(q-2)\frac{q^2}4$ planes that intersect $\pi_1$ and $\pi_2$ in skew lines.
 These include the $(q^2+1) - (2q+1)$ tangent planes of the points of $\mq \setminus (\pi_1 \cup \pi_2)$.
\end{itemize}

Now assume that $l \cap \mq$ are two distinct points $Q_1$ and $Q_2$.
Then $N_1, N_2 \notin l$.
We compute the numbers $p_{i,j}^2$.
\begin{itemize}
 \item $p^2_{1,1} = q-1$: For each point $P \in l$, the plane $\vspan{P,N_1,N_2}$ is the only plane through $P$ intersecting $\pi_1$ and $\pi_2$ in tangents.
 If $P \neq Q_1,Q_2$, then this plane intersects $\pi_1$ and $\pi_2$ in distinct points of $\mq$, hence is an oval plane.
 If $P$ equals $Q_1$ or $Q_2$, then this plane is the tangent plane to $P$.
 
 \item $p^2_{1,2} = (q-1)(\frac q 2 + 1)$: Take $P \in l$.
 Then $P$ lies on a unique tangent in $\pi_1$.
 The number of secants distinct from $l$ in $\pi_2$ through $P$ equals $q-1$ if $P \in \set{Q_1,Q_2}$, and $\frac q 2 - 1$ otherwise.
 \item $p^2_{1,3} = (q-1)(\frac q 2 - 1)$: Every point $P \in l$ lies on a unique tangent in $\pi_1$.
 If $P \neq Q_1, Q_2$, let $R$ be the unique point of $\mq$ on this tangent.
 Then $P$ lies on $\frac q 2$ skew lines in $\pi_2$.
 We disregard one of these, since it lies on the tangent plane through $R$.
 If $P \in \set{Q_1,Q_2}$, then $P$ lies on no skew lines.
 
 \item $p^2_{2,2} = \frac q 4 (q-1)(q+4)$:
 There are $q-1$ oval planes through $l$ distinct from $\pi_1$ and $\pi_2$.
 The number of secants distinct from $l$ through a point $P \in l$ in $\pi_1$ or $\pi_2$ equals $q-1$ if $P \in \set {Q_1, Q_2}$, and $\frac q 2 - 1$ otherwise.
 
 \item $p^2_{2,3} = \frac q 4 (q-1)(q-2)$: Take a point $P \in l$.
 If $P \in \set{Q_1,Q_2}$, then $P$ lies on no skew lines.
 Otherwise, $P$ lies on $\frac q 2 - 1$ secants distinct from $l$ in $\pi_1$, and $\frac q 2$ skew lines in $\pi_2$.
 
 \item $p^2_{3,3} = (q-1) \frac{(q-2)^2}4$: Every point $P \in l \setminus \set {Q_1,Q_2}$ lies on $\frac q 2$ skew lines in $\pi_1$ and $\pi_2$.
 Note that for each of the $(q^2+1)-2q$ points of $\mq \setminus (\pi_1 \cup \pi_2)$, we also counted its tangent plane.
\end{itemize}

Finally, assume that $l \cap \mq = \emptyset$.
We compute $p_{i,j}^3$.
Then no point of $l$ doesn't contain $N_1$, $N_2$, or any point of $\mq$.

\begin{itemize}
 \item $p^3_{1,1} = q+1$: Every point of $l$ lies on a unique tangent in $\pi_1$ and $\pi_2$.
 The plane through two such tangents intersects $\pi_1$ and $\pi_2$ in distinct points of $\mq$, hence is an oval plane.
 
 \item $p^3_{1,2} = (q+1)\frac q 2$: Every point of $l$ lies on a unique tangent in $\pi_1$ and $\frac q 2$ secants in $\pi_2$.
 \item $p^3_{1,3} = (q+1)\frac{q-4}2$: Take a point $P \in l$, and let $R$ be the unique point of the tangent $\vspan{P,N_1}$ in $\pi_1$.
 There are $\frac q 2 - 1$ skew lines distinct from $l$ through $P$ in $\pi_2$.
 One of these lines lies in the tangent plane through $R$.
 
 \item $p^3_{2,2} = (q+1) \frac{q^2}4$: Each point of $l$ lies on $\frac q 2$ secants in each of the planes $\pi_1$ and $\pi_2$.
 
 \item $p^3_{2,3} = (q+1) \frac q 4 (q-2)$: Each point of $l$ lies on $\frac q 2$ secants in $\pi_1$, and $\frac q 2 - 1$ skew lines distinct from $l$ in $\pi_2$.
 
 \item $p^3_{3,3} = \frac q 4 (q-4)(q-3) + 1$: There are $q-1$ other planes through $l$, and each point of $l$ lies on $\frac q 2 - 1$ skew lines distinct from $l$ in both of the planes $\pi_1$ and $\pi_2$.
 For each of the $(q^2 + 1) - 2(q+1)$ points of $\mq \setminus (\pi_1 \cup \pi_2)$, we have counted its tangent plane.
\end{itemize}

Thus, we have shown that the intersection numbers are well-defined, and therefore that $R_0, R_1, R_2, R_3$ constitutes an association scheme.
We can use the intersection numbers as described in \S \ref{SubsectionIntersectionMatrices} to construct the intersection matrices.
Consider the matrices
\begin{align*}
 P = \begin{pmatrix}
 1 & q^2-1 & \frac{q^2}2 (q+1) & \frac q 2 (q-1)(q-2) \\
 1 & q-1 & -q & 0 \\
 1 & -2 & q \frac{q-1}2 & -(q+1)\frac{q-2}2 \\
 1 & -(q+1) & 0 & q
 \end{pmatrix}, &&
 Q = \begin{pmatrix}
 1 & \frac q 2 (q^2+1) & q^2 & (q^2+1) \frac{q-2}2 \\
 1 & q \frac{q^2+1}{2(q+1)} & -2\frac{q^2}{q^2-1} & -(q-2)\frac{q^2+1}{2(q-1)} \\
 1 & - \frac{q^2+1}{q+1} & q\frac{q-1}{q+1} & 0 \\
 1 & 0 & -q\frac{q+1}{q-1} & \frac{q^2+1}{q-1} \\
 \end{pmatrix}.
\end{align*}
Then $Q$ diagonalises the intersection matrices, the corresponding eigenvalues are given in the columns of $P$, and $P Q = q(q^2+1) I_4$.
If we define the matrix $D = \text{diag}(1,q+1,q^2-1,q-1)$, and let $P_0, \dots, P_3$ denote the columns of $P$, then these assertions can be reformulated as $B_i Q D = Q D \text{diag}(P_i)$, and $P Q D = q(q^2+1) D$.
All these assertions are polynomial identities of degree at most 8.
Therefore, it suffices to check that they hold for at least 9 distinct values of $q$, which has been done by computer.

\subsection{Intersecting families in Möbius planes of even order}

Let $A_3$ denote the adjacency matrix corresponding to $R_3$.
Let $G$ denote the corresponding graph.
Then intersecting families in the Möbius plane are equivalent to cocliques in $G$.
We can find the eigenvalues of $A_3$ in the last column of $P$.
Hoffman's ratio bound tells us that
\[
 \alpha(G) \leq \frac{q^3 + q}{1 + \frac{\frac q 2 (q-1)(q-2)}{(q+1)\frac{q-2}2}  } = q(q+1).
\]
Note that each point of $\mq$ lies on one tangent plane and $q^2 + q$ oval planes.
Therefore, we know that this bound is tight, and we find the weak EKR property.
The $\tau$-eigenspace of $A_3$ is $V_2$, which has dimension $q^2$, as can be seen in the top row of $Q$.
Hence, by equality in Hoffman's ratio bound, if $K$ is an intersecting family of size $q(q+1)$ in the Möbius plane, then its characteristic vector lies in the $(q^2 + 1)$-dimensional space $\vspan{V_2,\one}$.

\bigskip

Let $W$ be the incidence matrix of the Möbius plane (see Definition \ref{DfIncidenceMatrix}).
Thus, the columns of $W$ are characteristic vectors of the sets of all circles through a point.
These are characteristic vectors of intersecting families of size $q(q+1)$, which implies that $\col(W) \leq \vspan{V_2,\one}$.
On the other hand, $W^t W$ is the square matrix indexed by the points of $\mq$, whose $(P,Q)$-entry gives the number of oval planes through $P$ and $Q$.
Since each plane through a secant line is an oval plane, this number equals $q+1$ if $P \neq Q$, and $q(q+1)$ if $P=Q$.
Thus, $W^t W = (q^2-1) I_{q^2+1} + (q+1) J_{q^2+1}$, which has full rank.
Since $\rk(W) = \rk(W^t W)$, this means that $\rk(W) = q^2 + 1 = \dim \vspan{V_2,\one}$.
Therefore, $\col(W) = \vspan{V_2,\one}$.

\bigskip

Hence, if $K$ is a maximum intersecting family, then its characteristic vector lies in $\col(W)$, or equivalently, equals $W \kappa$ for some vector $\kappa$.
Now take an oval plane $\pi \in K$, and let $W_\pi$ be the submatrix of $W$ obtained by deleting the rows corresponding to oval planes intersecting $\pi$ in at least one point of $\mq$.
Since $K$ only contains oval planes intersecting $\pi$, its characteristic vector has zeroes in all positions corresponding to oval planes not intersecting $\pi$, which index the rows of $W_\pi$.
This means that $W_\pi \kappa = 0$.
Note that every column in $W_\pi$ corresponding to a point of $\pi$ is a zero column.
Therefore, we also remove these columns from $W_\pi$, and let $\kappa_\pi$ denote the vector obtained by deleting the positions in $\kappa$ corresponding to the points of $\pi$.
Then it still holds that $W_\pi \kappa_\pi = 0$ since we only removed zero columns from $W_\pi$.

Consider the matrix $W_\pi^t W_\pi$.
It is a square matrix, indexed by the points of $\mq \setminus \pi$.
Its entry in the $(P,Q)$ position equals the number of oval planes through $P$ and $Q$ not intersecting $\pi$.
First suppose that $P = Q$.
There are $(q-1) \frac q 2$ skew lines in $\pi$, thus equally many planes through $P$ whose intersection with $\pi$ doesn't contain a point of $\mq$.
One of these planes is the tangent plane through $P$, thus the $(P,P)$-entry of $W_\pi^t W_\pi$ equals $(q-1) \frac q 2 - 1 = (q+1)\frac{q-2}2$.
Now suppose that $P \neq Q$.
The line $\vspan{P,Q}$ intersects $\pi$ in a point $R$ outside of $\mq$.
If $R$ is the nucleus $N$ of $\pi$, then every plane through $\vspan{P,Q}$ intersects $\pi$ in a tangent line.
If $R \neq N$, then there are $\frac q 2$ skew lines through $R$.
Note that no plane through the secant line $\vspan{P,Q}$ is tangent.
Therefore,
\[ W_\pi^t {W_\pi}_{P,Q} = \begin{cases}
 (q+1) \frac {q-2} 2 & \text{if } P = Q, \\
 0 & \text{if } P \neq Q \text{ and } \vspan{P,Q} \cap \pi = N, \\
 \frac q 2 & \text{otherwise.}
\end{cases} \]
The only tangent lines through $N$ lie in $\pi$.
This means that for every point $P \in \mq \setminus \pi$, $\vspan{P,N}$ is a secant line.
Thus, there is a unique point $Q \in \mq \setminus \pi$ such that $\vspan{P,Q} \cap \pi = N$.
We can choose in which way we order the columns of $W_\pi$.
This corresponds to ordering the points of $\mq \setminus \pi$.
Suppose that the ordering is $P_1, P_2, \dots, P_{q^2-q}$, such that $N \in \vspan{P_1,P_2}, N \in \vspan{P_3,P_4}, \dots, N \in \vspan{P_{q^2-q-1},P_{q^2-q}}$.
Then
\begin{align*}
 W_\pi^t W_\pi 
 & = I_\frac{q^2-q}2 \otimes \begin{pmatrix} (q+1)\frac{q-2}2 & 0 \\
 0 & (q+1)\frac{q-2}2 \\
 \end{pmatrix}
 + (J_\frac{q^2-q}2 - I_\frac{q^2-q}2) \otimes \left( \frac q 2 J_2 \right) \\
 & = I_\frac{q^2-q}2 \otimes \underbrace{\begin{pmatrix} (q+1)\frac{q-2}2 - \frac q 2 & - \frac q 2 \\
 - \frac q 2 & (q+1)\frac{q-2}2 - \frac q 2 \\
 \end{pmatrix}}_{=:M}
 + \frac q 2 J_{q^2-q}
\end{align*}

The eigenvalues of $M$ are given by $M \begin{pmatrix} 1 \\ 1 \end{pmatrix} = \left( (q+1)\frac{q-2}2 - q \right) \begin{pmatrix} 1 \\ 1 \end{pmatrix}$ and $M \begin{pmatrix} 1 \\ -1 \end{pmatrix} = (q+1)\frac{q-2}2 \begin{pmatrix} 1 \\ -1 \end{pmatrix}$.
Since, $q \geq 4$, the eigenvalues of $M$ are greater than zero.
As described in \S \ref{SubsectionKronecker}, $I \otimes M$ must then be positive definite.
Then $W_\pi^t W_\pi$ is the sum of a positive definite and a positive semi-definite matrix, and therefore itself positive definite.
In particular, the kernel of $W_\pi^t W_\pi$, which is equal to the kernel of $W_\pi$, is trivial.
Thus, $\kappa_\pi = 0$.
Hence, $\kappa(P) = 0$ for all $P \in \mq \setminus \pi$.

This means that for every point $P \in \mq$, if there is an oval plane $\pi \in K$ such that $P \notin \pi$, then the $P$-entry of $\kappa$ is zero.
Since $W \kappa \neq 0$, there must be a point $P$ on all oval planes of $K$.
Since $|K|$ equals the number of oval planes through $P$, $K$ must consist of all oval planes through $P$.
This proves Theorem \ref{ThmMainMobius}.

\bigskip

There is a unique Möbius plane of order 2 \cite[\S 3.2]{hartmann}.
This plane has 5 points, and 10 circles of size 3.
Therefore, all circles pairwise intersect, and no EKR property holds for this Möbius plane.

\subsection{2-intersecting families in Möbius planes}

Blokhuis and Bruen \cite{blokhuisbruen} proved that a 2-intersecting family in a Möbius plane of order $q$ has at most $\frac 1 2 q(q+1) + 1$ elements.
If we apply Hoffman's ratio bound to the weighted adjacency matrix $\frac{q+2}2 A_1 + A_3$ (see \cite[Theorem 2.4.2]{godsilmeagher} for details), then we also find this bound.
Unfortunately, the information provided by equality in Hoffman's ratio bound doesn't seem to tell us anything about the structure of a set attaining this bound, that cannot already be derived from the simple combinatorial proof by Blokhuis and Bruen. 

For ovoidal Möbius planes of odd order $q$, Blokhuis and Bruen improved the bound for 2-intersecting families to $\frac 1 2 q(q+1) - \frac q 4$.
They also proved that every ovoidal Möbius plane of order $q$ has a 2-intersecting family of size $3\frac{q-1}2 + 2$.
With a computer search, they found that in the classical Möbius planes (i.e.\ ovoidal Möbius planes constructed from an elliptic quadric) of order $q \leq 9$, 2-intersecting families are of size at most $2q$, with equality if and only if $q$ is even.

\section{Laguerre planes}
 \label{SectionLaguerre}

We characterise all intersecting families of maximum size in ovoidal Laguerre planes.
Recall that an ovoidal Laguerre plane of order $q$ is constructed from an oval cone in $\pg(3,q)$.
This cone has a unique singular point $R$, also called the vertex of the cone.
The oval planes are exactly the planes which don't contain $R$.

Again, we will prove that the relations from Definition \ref{DfRelAssoc} constitute an association scheme.
We need to make a distinction between Laguerre planes of odd and even order, because ovals behave quite differently in projective planes of odd and even order.
Before proceeding we make the following remarks which hold for all values of $q$.

\begin{rmk}
 \label{RmkIsomorphicOvalPlanes}
Let $\mq$ be an oval cone in $\pg(3,q)$ with vertex $R$.
Take two planes $\pi$ and $\rho$ not through $R$.
Consider the map $\Psi_{\pi,\rho}$ which maps a point $P \in \pi$ to the point $\vspan{R,P} \cap \rho$.
One easily checks that this is an isomorphism from $\pi$ to $\rho$, which maps $\mq \cap \pi$ to $\mq \cap \rho$.
Therefore, $\Psi_{\pi,\rho}$ maps an $i$-secant of $\mq \cap \pi$ to an $i$-secant of $\mq \cap \rho$.
Likewise, the number of $i$-secants to $\mq \cap \pi$ through $P \in \pi$ equals the number of $i$-secants to $\mq \cap \rho$ through $\Psi_{\pi,\rho}(P)$.
Also note that if $l$ is a line in $\pi$, then $\vspan{l, \Psi_{\pi,\rho}(l)}$ contains $R$, hence it is the unique non-oval plane through $l$.
\end{rmk}

\begin{rmk}
 \label{RmkOvoidalLaguerrePolynomials}
The classical construction of a Laguerre plane is as follows \cite[\S 4]{hartmann}:
Consider the set $U = \sett{f \in \FF_q[X]}{\deg f \leq 2}$ of polynomials over $\FF_q$ of degree at most 2.
For each polynomial $f(X) = aX^2 + bX + c \in U$, define $f(\infty) = a$.
Then we can see $f$ as a function $\FF_q \cup \set \infty \to \FF_q$.
Define the incidence structure $(\mathcal P, \mathcal B)$ where $\mathcal P = (\FF_q \cup \set \infty) \times \FF_q$ and the elements of $\mathcal B$ are the graphs $\sett{(x,f(x))}{x \in \FF_q \cup \set \infty}$ of all polynomials $f$ in $U$.
This incidence structure is isomorphic to the Laguerre plane arising from the quadric cone \cite[\S 4.4]{hartmann}.

It is not hard to check that both descriptions of the Laguerre plane are isomorphic.
There are coordinates $(X_0,X_1,X_2,X_3)$ of $\pg(3,q)$ such that the quadric cone $\mq$ has equation $X_0 X_2 = X_1^2$.
The vertex of $\mq$ is the point $R = (0,0,0,1)$.
The non-singular points of this quadric are of the form $(s^2,s,1,a)$ and $(1,0,0,a)$, with $s, a \in \FF_q$.
The oval planes are exactly the planes missing $R$, which are the planes with an equation of the form $X_3 = a X_0 + b X_1 + c X_2$.
We can identify the points $(s^2,s,1,a)$ and $(1,0,0,a)$ with $(s,a)$ and $(\infty,a)$ respectively in $(\FF_q \cup \set \infty) \times \FF_q$, and identify $X_3 = a X_0 + b X_1 + c X_2$ with (the graph of) the polynomial $a X^2 + b X + c$.
It's straightforward to check that this gives the desired isomorphism.

Thus, we have two descriptions of this Laguerre plane.
The \emph{geometrical description}, using the non-singular points and oval planes of a quadric cone, and the \emph{polynomial description}, using $(\FF_q \cup \set \infty) \times \FF_q$ and $U$.
\end{rmk}

\begin{rmk}
In a Laguerre plane of order $q$, the points are partitioned into $q+1$ parallel classes of $q$ points each.
Every circle contains a unique point of each parallel class.
In an ovoidal Laguerre plane, constructed from an oval cone $\mq$ with vertex $R$, the parallel classes are lines through $R$ with $R$ removed.
For each point $P$ we denote its parallel class as $\overline P$, which equals $\vspan{P,R} \setminus \set R$.
\end{rmk}

\subsection{The association scheme for ovoidal Laguerre planes of odd order}

Assume that $q$ is odd.
Since there is only one oval in $\pg(2,q)$ up to isomorphism, there is only one oval cone in $\pg(3,q)$ up to isomorphism, namely the quadric cone.
Thus, there is a unique ovoidal Laguerre plane of order $q$ and it has a geometrical and a polynomial description.
We will switch between the two descriptions, and use the most suitable one for our computations.

\bigskip

The first step is to prove that $R_0, R_1, R_2, R_3$ constitutes an association scheme.
In the polynomial interpretation, two distinct polynomials $f$ and $g$ are in relation $R_1$, $R_2$, or $R_3$ if and only if there are respectively 1, 2, or 0 values $x \in \FF_q \cup \set \infty$ such that $f(x)$ = $g(x)$.
(Strictly speaking the vertices of the association scheme are graphs of the polynomials, but we can identify each function with its graph.)

\bigskip

Take $a, b \in \FF_q^*$, $c \in \FF_q$, and $f_0 \in U$.
Consider the permutation
\[
 \sigma_{a,b,c,f_0}: U \to U: f(X) \mapsto a f(bX + c) + f_0(X).
\]
We also define the permutation
\[
 \varphi: U \to U: a X^2 + b X + c \mapsto c X^2 + b X + a.
\]
Let $\Sigma$ denote the group generated by the permutations $\sigma_{a,b,c,f_0}$ and $\varphi$ (with composition as operation).
Then $\Sigma$ acts on $U$.
In a natural way, $\Sigma$ also acts on $U \times U$, by defining $(f,g)^\sigma = (f^\sigma, g^\sigma)$.
The orbits of the action of $\Sigma$ on $U \times U$ are called the \emph{orbitals}.

\begin{lm}
 \label{LmOrbitalsSigma}
 The orbitals of $\Sigma$ are the relations $R_0, R_1, R_2, R_3$.
\end{lm}

\begin{proof}
Take two polynomials $f$ and $g$ in $U$.
Take $\sigma = \sigma_{a,b,c,f_0}$ in $\Sigma$ and $x \in \FF_q$.
Then $f^\sigma(x) = g^\sigma(x)$ if and only if $f(bx+c) = g(bx+c)$, and $f^\sigma(\infty) = g^\sigma(\infty)$ if and only if $f(\infty) = g(\infty)$.
Since $X \mapsto b X + c$ is a permutation of $\FF_q$, we see that $(f,g)^\sigma$ is in the same relation $R_i$ as $(f,g)$.
Analogously, $f^\phi(x) = g^\phi(x)$ if and only if $f(x^{-1}) = g(x^{-1})$ (with the convention that $0$ and $\infty$ are each others inverses).
Therefore, the relations $R_i$ are unions of orbitals of $\Sigma$.

Now consider the subgroup $\Sigma_0$ of $\Sigma$, consisting of the elements $\sigma_{a,b,c,0}$.
Partition $U$ in the sets $Q_{-1} = \set 0$, $Q_0$ are non-zero the constant polynomials, $Q_1$ are the polynomials of degree 1, and for $i=0,1,2$, $Q_2^{i}$ are the polynomials of degree 2, with $i$ distinct roots.
We claim that these 6 sets are the orbits of the action of $\Sigma_0$ on $U$.
For every polynomial $f$ and every $\sigma \in \Sigma_0$, $\deg f^\sigma = \deg f$.
This means that $0^\sigma = 0$, and therefore $(f,0)$ and $(f,0)^\sigma = (f^\sigma,0)$ are in the same relation $R_i$.
In other words, every element of $\sigma$ leaves the degree and the number of roots over $\FF_q \cup \set \infty$ invariant.
Thus, $\Sigma_0$ leaves the 6 set in which we partitioned $U$ invariant, which means they are unions of orbits of $\Sigma_0$.
Therefore, it suffices to prove that $\Sigma_0$ works transitively on each of these sets.

It is obvious that $Q_{-1}$ and $Q_0$ are orbits of $\Sigma_0$.
Each polynomial $f \in Q_1$ is of the from $b X + c = (X)^{\sigma_{1,b,c,0}}$.
Each polynomial $f \in Q_2^1$ is of the form $a (X+c)^2 = (X^2)^{\sigma_{a,1,c,0}}$.
Thus, $\Sigma_0$ works transitively on $Q_1$ and $Q_2^1$.
Now take a non-zero element $\gamma \in \FF_q^*$ and suppose that $(X^2-\gamma)^{\sigma_{a,b,c,0}} = X^2 - \gamma$.
Then $ab^2 X^2 + 2abc X + a(c^2 - \gamma) = X^2 - \gamma$, or equivalently
\[ \begin{cases}
 ab^2 = 1, \\
 2abc = 0, \\
 a(c^2-\gamma) = -\gamma.
\end{cases}\]
This is equivalent to $c=0$, $a=1$, and $b=\pm 1$, hence there are two choices.
Therefore, the orbit of $X^2 - \gamma$ has size $\frac{|\Sigma_0|}2 = q\frac{(q-1)^2}2$.
On the other hand, $a X^2 + b X + c$ is in $Q_2^0$ (respectively $Q_2^2)$ if and only if $a$ is non-zero and $D = b^2 - 4 a c$ is not a square (respectively a non-zero square).
From this we can calculate that $|Q_2^0| = |Q_2^2| = q \frac{(q-1)^2}2$.
$X^2 - \gamma$ is in $Q_2^0$ or $Q_2^2$ depending on whether $\gamma$ is a square.
This proves that $\Sigma_0$ acts transitively on $Q_2^0$ and $Q_2^2$.

Now let $\Sigma_1$ denote the stabiliser of the zero polynomial in $\Sigma$.
Then $\Sigma_1$ contains $\Sigma_0$ and $\varphi$.
Note that $\varphi$ maps a constant non-zero polynomial to an element of $Q_2^1$, and maps $X+1$ to $X(X+1) \in Q_2^2$.
Therefore, the orbits of the action of $\Sigma_1$ on $U$ are unions of $Q_{-1}$, $Q_2^0$, $Q_0 \cup Q_2^1$, and $Q_1 \cup Q_2^2$.
By previous considerations $f$ and $f^\sigma$ are in the same relation w.r.t.\ 0 for each $\sigma \in \Sigma_1$.
Therefore $Q_{-1}$, $Q_2^0$, $Q_0 \cup Q_2^1$, and $Q_1 \cup Q_2^2$ are exactly the orbits of $\Sigma_1$.

To finish the proof, we need to show that if $(f,g)$ and $(h,p)$ are in the same relation, then $(f,g)$ and $(h,p)$ are in the same orbital of $\Sigma$.
Note that $(f,g)$ is in the same orbital as $(f-g,0)$ (apply $\sigma_{1,1,0,-g}$), and likewise $(h,p)$ is in the same orbital as $(h-p,0)$.
This means that $f-g$ and $h-p$ are in the same relation w.r.t.\ 0, which means that $f-g = (h-p)^\sigma$ for some $\sigma \in \Sigma_1$.
Therefore $(f-g,0)$ and $(h-p,0)$ are in the same orbital of $\Sigma$.
\end{proof}

Suppose that $(f,g)$ and $(h,p)$ are both in $R_k$.
Then $(h,p) = (f,g)^\sigma$.
A polynomial $t$ satisfies $(f,t) \in R_i$ and $(t,g) \in R_j$ if and only if $(h,t^\sigma) \in R_i$ and $(t^\sigma,p) \in R_j$.
This shows that the intersection numbers are well-defined, and therefore proves that $R_0, R_1, R_2, R_3$ is an association scheme.
To find the matrices $P$ and $Q$, we still need to compute the intersection numbers.
But knowing that the relations $R_i$ are an association scheme allows us to use the known identities.
Let $n_i$ denote the number $p^0_{i,i}$, called the \emph{valency} of $R_i$.
Then $n_1 = |Q_0 \cup Q_2^1| = q^2-1$, $n_2 = |Q_1 \cup Q_2^2| = q\frac{q^2-1}2$, and $n_3 = |Q_2^0| = q \frac{(q-1)^2}2$

\begin{lm}[{\cite[Lemma 2.1.1]{bcn}}]
 \label{LmIntersectionIdentities}
 \begin{enumerate}
  \item $p_{0,j}^k = \delta_{j,k}$,
  \item $p_{i,j}^k n_k = p_{i,k}^j n_j$,
  \item $\displaystyle \sum_{j=0}^3 p_{i,j}^k = n_i$.
 \end{enumerate}
\end{lm}

We determine the intersection numbers in the geometrical description of the ovoidal Laguerre plane.
Take two oval planes $\pi$ and $\rho$, and let the line $l$ be their intersection.
Let $\Psi$ denote the map $\Psi_{\pi,\rho}$ defined in Remark \ref{RmkIsomorphicOvalPlanes}.
Note that $\Psi(l) = l$.

\begin{itemize}
 \item $p^1_{1,1} = q-2$: Assume that $l$ contains a unique point $Q$ of $\mq$.
 There are $q-2$ oval planes through $l$, distinct from $\pi$ and $\rho$, which are all in relation $R_1$ to $\pi$ and $\rho$.
 Suppose that $\tau$ is a plane not through $l$ that intersects $\pi$ and $\rho$ in a unique point of $\mq$.
 Then $\tau$ intersects $l$ in a point $P$, which must lie on a tangent line $l_1 \neq l$ in $\pi$ and $l_2 \neq l$ in $\rho$.
 Then $P \neq Q$, so $P$ is external and $l_i$ is the only tangent distinct from $l$ through $P$ in $\pi$. 
 Thus $\Psi(l_1)$ must be $l_2$, which implies that $\tau$ contains $R$.
 Hence, $p_{1,1}^1 = q-2$.
 
 \item $p^2_{1,1} = q-1$: Now assume that $l$ intersect $\mq$ in two points $Q_1$ and $Q_2$.
 Then $Q_1$ only lies on unique tangents $l_1 \neq l$ in $\pi$ and $l_2 \neq l$ in $\rho$.
 As above, the plane $\vspan{l_1,l_2}$ contains $R$ so $Q_1$ (and likewise $Q_2$) lies on no oval planes intersecting $\pi$ or $\rho$ in a tangent.
 There are $\frac{q-1}2$ external points on $l$ (by applying $\Psi$ we see that a point of $l$ is external to $\mq \cap \pi$ if and only if it is external to $\mq \cap \rho$).
 Take such a point $P$.
 There are two choices for a tangent line $t$ through $P$ in $\pi$.
 Then $P$ lies on a unique tangent line distinct from $\Psi(t)$ in $\rho$.
 Thus, each tangent line through $P$ gives a unique oval plane intersecting $\pi$ and $\rho$ in a tangent.
 Therefore, $p^2_{1,1} = 2 \frac{q-1}2$.
 
 \item $p^2_{1,2} = \frac{q^2-1}2$: Assume again that $\mq \cap l = \set{Q_1,Q_2}$.
 We count the number of planes $\tau$ which intersect $\pi$ in a tangent line $l_1$ and $\rho$ in a secant line $l_2$.
 Then $\Psi(l_1)$ cannot be $l_2$, and $\tau$ must be an oval plane.
 Note that $\tau$ cannot contain $l$.
 Therefore, it suffices to take the sum over all points $P$ on $l$ of the number of tangent lines $l_1$ through $P$ in $\pi$ times the number of secant lines $l_2 \neq l$ through $P$ in $\rho$. 
 
 Thus, take a point $P$ on $l$.
 First assume that $P \in \set{Q_1,Q_2}$.
 Then $P$ lies on a unique tangent in $\pi_1$, and $q-1$ secants $l_2 \neq l$ in $\rho$.
 Now assume that $P$ is one of the $\frac{q-1}2$ external points.
 Then $P$ lies on 2 tangents in $\pi$ and $\frac{q-1}2 - 1$ secants $l_2 \neq l$ in $\rho$.
 We can disregard internal points, since they do not lie on tangents.
 Thus,
 \[
  p^2_{1,2} = 2 \cdot 1 \cdot (q-1) + \frac{q-1}2 \cdot 2 \cdot \left( \frac{q-1}2 - 1 \right) = \frac{q^2-1}2.
 \]
 
 \item $p^3_{2,3} = (q+1)\frac{(q-1)^2}4$:
 Assume that $l$ is skew to $\mq$.
 Similarly as the previous calculation, we need to take the sum over all points $P \in l$ of the number of secants through $P$ in $\pi$ multiplied with the number of skew lines distinct from $l$ through $P$ in $\rho$.
 $l$ contains $\frac{q+1}2$ external points, which lie on $\frac{q-1}2$ secants in $\pi$ and $\frac{q-1}2$ skew lines in $\rho$, and $\frac{q+1}2$ internal points which lie on $\frac{q+1}2$ secants in $\pi$ and $\frac{q+1}2$ skew lines in $\rho$.
 Thus,
 \[
  p^3_{2,3} = \frac{q+1}2 \cdot \frac{q-1}2 \left( \frac{q-1}2 - 1 \right)
  + \frac{q+1}2 \cdot \frac{q+1}2 \left( \frac{q+1}2 - 1 \right)
  = (q+1) \frac{(q-1)^2}4.
 \]
\end{itemize}

These four intersection numbers, together with the identities from Lemma \ref{LmIntersectionIdentities} suffice to calculate all intersection numbers, or equivalently the intersection matrices.

Consider the matrices
\begin{align*}
P = Q = \begin{pmatrix}
1 & q^2-1 & q \frac{q^2-1}2 & q\frac{(q-1)^2}2 \\
1 & -1 & q \frac{q-1}2 & -q \frac{q-1}2 \\
1 & q-1 & -q & 0\\
1 & -(q+1) & 0 & q
\end{pmatrix}.
\end{align*}
Then $Q$ diagonalises the intersection matrices, the columns of $P$ give the corresponding eigenvalues, and $PQ = q^3I_4$.
Therefore, $P$ and $Q$ are the matrices of eigenvalues and dual eigenvalues of the association scheme.
These assertions have been checked by formulating them as polynomial identities of degree at most 6 and checking that they hold for at least 7 values of $q$.

The association scheme in question is an example of a \emph{self-dual translation scheme}, see \cite[\S 2.10B]{bcn}.
Self-dual refers to $P=Q$, translation scheme refers to the fact that the graphs in the association scheme are Cayley graphs on an abelian group.

\subsection{Intersecting families in ovoidal Laguerre planes of odd order}

Let $A_i$ denote the adjacency matrices of the association scheme, and $V_i$ the eigenspaces.
Let $G$ denote the graph with adjacency matrix $A_3$.
The cocliques of $G$ are exactly the intersecting families of circles in the ovoidal Laguerre plane.
We apply Hoffman's ratio bound (the eigenvalues of $A_3$ can be found in the last column of $P$).
\[
 \alpha(G) \leq \frac{q^3}{ 1 + \frac{ q\frac{(q-1)^2}2 }{ q \frac{q-1}2 }} = q^2.
\]
The matrix $A_3$ attains its smallest eigenvalue on the eigenspace $V_1$.
This space has dimension $q^2-1$, as can be read from the top row of $Q$.

The number of planes through a point $P$ in $\pg(3,q)$ that miss a fixed point $R$ also equals $q^2$.
Thus, we have proven the weak EKR property.
Now we prove the strong EKR property, using the polynomial description.

\bigskip

Let $W$ denote the incidence matrix of the Laguerre plane (see Definition \ref{DfIncidenceMatrix}).
Consider the matrix $W W^t$.
Its rows and columns are indexed by the polynomials in $U$, and give the number of points on which the polynomials take the same value.
Therefore,
\[
 W W^t = (q+1) I + A_1 + 2 A_2.
\]
Its eigenvalues are given by $P (q+1, 1, 2, 0)^t$.
One can calculate from this that $\ker(W W^t) = \vspan{V_2,V_3}$, and $\col(W W^t) = \col(W) = \vspan{\one, V_1}$.

Now take an intersecting family of size $q^2$ in the Laguerre plane.
By equality in Hoffman's ratio bound, its characteristic vector lies in $\vspan{\one,V_1} = \col(W)$.
Thus, this characteristic vector equals $W \kappa$ for some vector $\kappa$, whose entries are indexed by $(\FF_q \cup \set \infty) \times \FF_q$.

Take a polynomial $f \in K$.
Let $W_f$ denote the restriction of $W$ to the rows corresponding to polynomials $g$ that don't intersect $f$, or equivalently $g-f \in Q_2^0$.
Since these polynomials $g$ cannot be in $K$, $W_f \kappa = 0$.
Each column of $W$ corresponding to a point $(x,f(x))$ is a zero column.
Therefore, if we also remove these columns from $W_f$ and the let $\kappa_f$ denote the restriction of $\kappa$ to points $(x,y)$ with $y \neq f(x)$, then $W_f \kappa_f  = 0$.
The next step is to determine $\ker(W_f) = \ker(W_f^t W_f)$.
The $((x_1,y_1), (x_2,y_2))$ entry of $W_f^t W_f$ equals the number of $g \in U$ such that $g-f \in Q_2^0$, $g(x_1) = y_1$, and $g(x_2) = y_2$.
By considering $h = g-f$, this number equals the number of $h \in Q_2^0$ such that $h(x_1) = y_1 - f(x_1)$, and $h(x_2) = y_2 - f(x_2)$.
Note that $y_1 - f(x_1)$ and $y_2 - f(x_2)$ are non-zero.

\begin{lm}
 \label{LmQ20ThroughPoint}
 Take $x \in \FF_q \cup \set \infty$, and $y \in \FF_q^*$.
 Then the number of $h \in Q_2^0$ with $h(x) = y$ equals $q \frac{q-1}2$.
\end{lm}

\begin{proof}
Take $h \in Q_2^0$.
Then there is a unique $a \in \FF_q^*$ such that $ah(x) = y$.
Therefore, the number of $h \in Q_2^0$ with $h(x) = y$ equals $|Q_2^0|/(q-1) = q \frac{q-1}2$.
\end{proof}

Now take two points $(x_1,y_1)$ and $(x_2,y_2)$ in $(\FF_q \cup \set \infty) \times \FF_q^*$, with $x_1 \neq x_2$.
Let $m_i$ denote the number of $h \in U$ such that $h(x_1) = y_1$, $h(x_2) = y_2$, and $h$ has exactly $i$ roots over $\FF_q \cup \set \infty$.
Let $S$ and $\overline S$ denote the sets of non-zero squares and of non-squares of $\FF_q$ respectively.

\begin{lm}
 \label{LmQOddQ1}
\begin{enumerate}
 \item[(1)] $ \displaystyle m_1 = \begin{cases}
 2 & \text{if } \frac{y_1}{y_2} \in S, \\
 0 & \text{if } \frac{y_1}{y_2} \in \overline S,
 \end{cases} $
\item[(2)] $m_0 + m_1 + m_2 = q$,
\item[(3)] $m_1 + 2 m_2 = q-1$.
\end{enumerate}
\end{lm}

\begin{proof}
(1) The number $m_1$ gives the number of $h \in Q_0 \cup Q_2^1$ such that $h(x_1) = y_1$, and $h(x_2) = y_2$.
There is one such $h \in Q_0$ if $x_1, x_2 \neq \infty$ and $y_1 = y_2$, and no such $h$ otherwise.

Next we determine the number of $h \in Q_2^1$.
These are of the form $h(X) = a(X-b)^2$, with $a \neq 0$.

First suppose that $x_1, x_2 \neq \infty$.
We are interested in the number of solutions to the system of equations
\[
 \begin{cases}
 y_1 = a(x_1-b)^2, \\
 y_2 = a(x_2-b)^2.
 \end{cases}
\]
We can rewrite this as
\[ \begin{cases}
 \frac{y_1}{y_2} = \left( \frac{x_1-b}{x_2-b} \right)^2, \\
 a = \frac{y_1}{(x_1-b)^2}.
\end{cases} \]
If we find a solution $b$ for the first equation, then there exists a unique $a$ that solves the second equation.
Hence, we are interested in the number of solutions to the first equation.
Evidently, there are no solutions if $y_1/y_2$ is not a square.
If $y_1=y_2$, then this equation reduces to $(x_1-b)^2=(x_2-b)^2$.
Since $x_1\neq x_2$, this is equivalent to $x_1-b = b-x_2$, which has a unique solution.
Therefore, the only remaining case is the case where $y_1/y_2$ is a square $\alpha^2 \neq 0, 1$.
Note that $\frac{x_1-b}{x_2-b} = \frac{x_1-x_2}{x_2-b} + 1$.
Therefore, we are looking for values of $b$ such that $\frac{x_1-x_2}{x_2-b} + 1 = \pm \alpha$, or equivalently
$b = x_2-\frac{x_1-x_2}{\pm \alpha - 1}$.
This permitted because $\alpha^2 \neq 1$, so $\pm \alpha -1 \neq 0$.
This way we see that there are 2 solutions for $b$.

Secondly suppose that $x_1$ or $x_2$ equals $\infty$.
Without loss of generality, we may assume that $x_1 = \infty$.
Then $m_1$ is the number of polynomials $h = y_1(X-b)^2 \in Q_2^1$ such that $y_1(x_2-b)^2 = y_2$, or equivalently $(x_2-b)^2 = y_2/y_1$.
There are two options if $y_2/y_1$ is a square, and none otherwise.

\bigskip

Now take a value $x \neq x_1, x_2$ in $\FF_q \cup \set \infty$.
Consider the Vandermonde matrix
\[
 V = \begin{pmatrix}
 x^2 & x & 1 \\
 x_1^2 & x_1 & 1 \\
 x_2^2 & x_2 & 1
 \end{pmatrix}.
\]
If $x$, $x_1$, or $x_2$ equals infinity, then replace the row $\begin{pmatrix} \infty^2 & \infty & 1 \end{pmatrix}$ by $ \begin{pmatrix} 1 & 0 & 0 \end{pmatrix}$.
The matrix $V$ is non-singular.

\bigskip

(2) $m_0 + m_1 + m_2$ gives the number of polynomials $h = a X^2 + b X + c$ such that $h(x_1) = y_1$ and $h(x_2) = y_2$.
Take $x \in (\FF_q \cup \set \infty) \setminus \set{x_1,x_2}$.
For every value of $y \in \FF_q$, there is a unique such $h$ with $h(x) = y$, since $V (
a,b,c )^t = (y,y_1,y_2)^t$ has a unique solution for $a$, $b$, and $c$.
Thus, the number in question equals $q$.

\bigskip

(3) Consider the set 
\[
Z = \sett{(h,x)}{h \in U, \, x \in \FF_q \cup \set \infty, \, h(x) = 0, \, h(x_1) = y_1, \,  h(x_2) = y_2}.
\]
Given a value $x$, the number of solutions to $V (a, b, c )^t = ( 0 , y_1 , y_2)^t$ is 1 if $x \neq x_1, x_2$, and 0 otherwise.
Therefore, $|Z| = q-1$.
On the other hand, $|Z|$ is clearly $m_1 + 2m_2$.
\end{proof}

\begin{table}[ht!]
    \centering
    \begin{tabular}{c|c c c}
         &  $\overline S$ & $S$ \\ \hline
         $m_0$ & $\frac{q+1}2$ & $\frac{q-1}2$ \\
         $m_1$ & 0&2 \\
         $m_2$ & $\frac{q-1}2$ & $\frac{q-3}2$
    \end{tabular}
    \caption{Given two points $(x_1,y_1)$ and $(x_2,y_2)$ in $(\FF_q \cup \set \infty) \times \FF_q^*$ with $x_1 \neq x_2$, this table gives the number $m_i$ of polynomials $h$ with $i$ roots over $\FF_q \cup \set \infty$ such that $h(x_1) = y_1$ and $h(x_2) = y_2$, depending on wether $y_1/y_2 \in S$ or $y_1/y_2 \in \overline S$.}
    \label{TableMi}
\end{table}

This gives us Table \ref{TableMi}.
This implies that
\[
 W_f^t W_f((x_1,y_1), (x_2,y_2)) = \begin{cases}
 q \frac{q-1}2 & \text{if } (x_1,y_1) = (x_1,y_2), \\
 0 & \text{if } x_1 = x_2 \text{ and } y_1 \neq y_2, \\
 \frac{q-1}2 & \text{if } x_1 \neq x_2 \text{ and } \frac{y_1 - f(x_1)}{y_2 - f(x_2)} \in S, \\
 \frac{q+1}2 & \text{if } x_1 \neq x_2 \text{ and } \frac{y_1 - f(x_1)}{y_2 - f(x_2)} \in \overline S.
 \end{cases}
\]

Now suppose that the ordering of the columns of $W_f$ is as follows: there are $q+1$ consecutive blocks consisting of the $q-1$ points with the same $x$-coordinate.
In each such block, the first $\frac{q-2}2$ points $(x,y)$ satisfy $y-f(x) \in S$, the next $\frac{q-1}2$ points satisfy $y-f(x) \in \overline S$.
Then
\begin{align*}
 W_f^t W_f & = q\frac{q-1}2 I_{q^2-1} + (J_{q+1} - I_{q+1}) \otimes \underbrace{\left(\begin{pmatrix} \frac{q-1}2 & \frac{q+1}2 \\ \frac{q+1}2 & \frac{q-1}2 \end{pmatrix} \otimes J_{\frac{q-1}2} \right)}_{ =: M }.
\end{align*}

We use \S \ref{SubsectionKronecker} to determine $\ker(W_f^t W_f)$.
This kernel equals the eigenspace of $(J-I) \otimes M$ with eigenvalue $- q \frac{q-1}2$.
\begin{align*}
 \begin{pmatrix} \frac{q-1}2 & \frac{q+1}2 \\ \frac{q+1}2 & \frac{q-1}2 \end{pmatrix}
 \begin{pmatrix} 1 \\ 1 \end{pmatrix}
 = q \begin{pmatrix} 1 \\ 1 \end{pmatrix},
 && \begin{pmatrix} \frac{q-1}2 & \frac{q+1}2 \\ \frac{q+1}2 & \frac{q-1}2 \end{pmatrix}
 \begin{pmatrix} 1 \\ -1 \end{pmatrix}
 = - \begin{pmatrix} 1 \\ -1 \end{pmatrix}.
\end{align*}
The eigenvectors of $J_\frac{q-1}2$ are $\one$ with eigenvalue $\frac{q-1}2$ and the vectors of $\vspan{\one}^\perp$ with eigenvalue 0.
Therefore, $M$ has the following eigenvectors with non-zero eigenvalue:
\begin{itemize}
 \item $\one_{q-1}$ with eigenvalue $q\frac{q-1}2$,
 \item $\begin{pmatrix} 1 \\ -1 \end{pmatrix} \otimes \one_\frac{q-1}2$ with eigenvalue $-\frac{q-1}2$.
\end{itemize}
$J_{q+1} - I_{q+1}$ has eigenspaces $\vspan{\one_{q+1}}$ with eigenvalue $q$ and $\vspan{\one_{q+1}}^\perp$ with eigenvalue $-1$.
Let $v_1, \dots, v_q$ be a basis for $\vspan{\one_{q+1}}^\perp$.
Then $v_1 \otimes \one_{q-1}, \dots, v_q \otimes \one_{q-1}$ and $\one_{q+1} \otimes \begin{pmatrix} 1 \\ -1 \end{pmatrix}  \otimes \one_\frac{q-1}2$ are a basis for the eigenspace of $(J-I) \otimes M$ with eigenvalue $-q\frac{q-1}2$, which equals $\ker(W_f^t W_f) =\ker(W_f)$.
Thus, 
$$\kappa_f = a_f \one_{q+1} \otimes \begin{pmatrix} \one_\frac{q-1}2 \\ -\one_\frac{q-1}2 \end{pmatrix} + v_f \otimes \one_{q-1},$$
where $a_f$ is a constant and $v_f \in \vspan{\one_{q+1}}^\perp$.
This holds for every $f \in K$.
We can interpret this as follows: $a_f$ is a constant and $v_f$ is a vector, whose positions are indexed by $\FF_q \cup \set \infty$ and whose entries sum to zero, such that $\kappa(x,y)$ equals $a_f + v_f(x)$ if $y - f(x) \in S$ and $-a_f + v_f(x)$ if $y-f(x) \in \overline S$.
We don't get any information about $\kappa(x,f(x))$.

\begin{lm}
 \label{LmDifferentSquareTypes}
If $q$ is odd, then for every two distinct elements $a$ and $b$ in $\FF_q$, there exists an element $c$ such that of the two numbers $c-a$ and $c-b$, exactly one belongs to $S$ and the other one to $\overline S$.
\end{lm}

\begin{proof}
First assume that $q \equiv 1 \pmod 4$.
Then the number $c$ must be a neighbour of exactly one of the elements $a$ and $b$ in the Paley graph of order $q$ (see e.g.\ \cite[\S 5.8]{godsilmeagher}).
There are $\frac{q-5}4$ common neighbours of $a$ and $b$, and $\frac{q-1}4$ common non-neighbours of $a$ and $b$, which leaves $q - 2 - \frac{q-5}4 - \frac{q-1}4 = \frac{q-1}2$ elements which are a neighbour of $a$ or $b$ but not both.

Now assume that $q \equiv 3 \pmod 4$.
Then $c-a$ and $c-b$ are of different quadratic type if and only if $z := \frac{c-a}{b-a}$ and $\frac{c-b}{b-a} = z-1$ are of different quadratic type.
Thus, we need to find a number $z$ such that $z$ and $z-1$ are non-zero and of different quadratic type.
Let $p$ denote the characteristic of $\FF_q$.
Consider the pairs $(1,2), (2,3), \dots, (p-2,p-1)$ in $\FF_q$.
Since $q \equiv 3 \pmod 4$, $1$ is a square and $-1$ is a non-square.
So at least one of these pairs consists of a square and a non-square, which give us the desired numbers $z-1$ and $z$.
\end{proof}

\begin{lm}
For every $f \in K$, $a_f = 0$.
\end{lm}

\begin{proof}
Take $f$ and $g$ in $K$.
Then there exists some $x \in \FF_q \cup \set \infty$ for which $f(x) = g(x)$.
Choose $y_1$ and $y_2$ such that $y_1 - f(x) \in S$ and $y_2 - f(x) \in \overline S$.
Then
\begin{align*}
 \kappa(x,y_1) &= a_f + v_f(x) = a_g + v_g(x), \\
 \kappa(x,y_2) &= - a_f + v_f(x) = -a_g + v_g(x).
\end{align*}
Therefore, $\kappa(x,y_1) - \kappa(x,y_2) = 2 a_f = 2 a_g$.
This means that $a_f$ is equal for all $f \in K$.
This also means that $v_f(x) = v_g(x)$ if $f(x) = g(x)$.
Suppose that $a_f = a$ for all $f \in K$.

There exist $f$ and $g$ in $K$ which are intersecting, but not 2-intersecting.
Otherwise, $K$ would be a 2-intersecting family, and $|K| < q^2$ (see Theorem \ref{ThmLaguerre2intersecting}).
There are $q$ values of $x \in \FF_q \cup \set \infty$ for which $f(x) \neq g(x)$.
Take such a value $x$, then by Lemma \ref{LmDifferentSquareTypes} there exists a $y \in \FF_q$ such that of the two numbers $y-f(x)$ and $y-g(x)$, one is in $S$, and the other one is in $\overline S$.
Thus, $\kappa(x,y) = v_f(x) + \varepsilon a = v_g(x) - \varepsilon a$ with $\varepsilon = \pm 1$.
This means that $v_f(x) - v_g(x) = \pm 2a$.
Let $n_+$ denote the number of $x$ for which $v_f(x) - v_g(x) = 2a$.
Then there are $q-n_+$ values of $x$ for which $v_f(x) - v_g(x) = -2a$, and one for which $v_f(x) = v_g(x)$.
Then $(v_f - v_g) \cdot \one_{q+1} = n_+ \cdot 2a + (q-n_+) \cdot (-2a) = 2a (2 n_+ - q)$.
On the other hand, $v_f$ and $v_g$ both lie in $\vspan{\one}^\perp$, which means that $2a (2 n_+ - q) = (v_f - v_g) \cdot \one = 0$.
If $a \neq 0$, then $n_+ = q/2$, which is impossible since $q$ is odd and $n_+$ is an integer.
Hence, $a=0$.
\end{proof}

The previous lemma implies that $\kappa_f = v_f \otimes \one_{q-1}$ for every $f \in K$, $v_f \in \vspan{\one}^\perp$.
This means that for every $x \in \FF_q \cup \set \infty$ and every $f \in K$, $\kappa$ takes the value $v_f(x)$ on all points $(x,y)$ with $y \neq f(x)$ and $\sum_x v_f(x) = 0$.

Take two functions $f$ and $g$ in $K$.
For each $x \in \FF_q \cup \set \infty$ there exists a $y$ with $f(x) \neq y \neq g(x)$, since $q \geq 3$.
Then $v_f(x) = \kappa(x,y) = v_g(x)$.
This means that $v_f$ is equal for all $f \in K$, and we drop the subscript.
Assume that for each point $(x,y) \in (\FF_q \cup \set \infty) \times \FF_q$ there exists a polynomial $f \in K$ with $f(x) \neq y$.
Then for each point $(x,y)$, $\kappa(x,y) = v(x)$.
Now take the row $w$ in $W$ corresponding to a certain polynomial $g$.
Then $w \kappa = \sum_{x \in \FF_q \cup \set \infty} \kappa(x,g(x)) = \sum_x v(x) = 0$, since $v \in \vspan{\one}^\perp$.
But then $w \kappa = 0$ for each row $w$ of $W$ which means that the characteristic vector of $K$ is $W \kappa = 0$.
This is a contradiction, so there must exist a point $(x,y)$ such that $f(x) = y$ for all $f \in K$.
Since there are only $q^2$ polynomials that satisfy this condition, $K$ consists of all polynomials $f$ with $f(x) = y$.
This proves Theorem \ref{ThmMainLaguerreOdd}.

\subsection{2-intersecting families in Laguerre planes}
 \label{SubsectionLaguerre2intersecting}

The spectral methods to bound the size of a 2-intersecting family in a Laguerre plane of odd order (the Delsarte LP clique and coclique bounds) yield as bound $\frac{q^2+1}2$.
We can also prove this bound using more or less the same arguments as Blokhuis and Bruen \cite{blokhuisbruen} used for Möbius planes.
The advantage of the second method is that this works for all Laguerre planes, only uses elementary arguments, and can easily be proven to not be feasible.

\bigskip

So suppose that $K$ is a 2-intersecting family in a Laguerre plane of order $q$.
Take a point $P$ and look in the residue at $P$.
The circles through $P$ form an intersecting family in this residue.
This residue is an affine plane of order $q$, and the circles through $P$ are the lines in this affine plane, excluding one parallel class.
Since lines intersect if and only if they are in a different parallel class, and there are $q+1$ parallel classes in an affine plane of order $q$, this means that there are at most $q$ circles through $P$ in $K$.
Now fix a circle $c \in K$.
Perform a double count on the set $Z = \sett{(P,c')}{c' \in K \setminus \set c, \, P \in c \cap c'}$.
Since $K$ is 2-intersecting, $|Z| = 2(|K| - 1)$.
On the other hand, $c$ contains $q+1$ points, and every point lies on at most $q-1$ other circle of $K$, thus $|Z| \leq (q+1)(q-1)$.
Therefore, $|K| \leq 1 + \frac{(q+1)(q-1)}2 = \frac{q^2+1}2$.

If this bound were tight, then every point would lie on either $q$ or 0 circles of $K$.
Therefore, the number of points covered by the circles of $K$ would equal
\[
 \frac{(q+1) \cdot \frac{q^2+1}2}q = \frac 1 2 \left( q^2 +q + 1 + \frac 1 q \right),
\]
which isn't an integer.
Thus we have proven the following theorem

\begin{thm}
 \label{ThmLaguerre2intersecting}
If $\mathcal F$ is a 2-intersecting family in a Laguerre plane of order $q$, then $|\mathcal F| < \frac{q^2+1}2$.
\end{thm}

For ovoidal Laguerre planes of even order, this bound can be improved substantially.
In section \S \ref{SubsectionLaguerre2intersectingEven} we prove a weak EKR bound for these Laguerre planes.
For ovoidal Laguerre planes of order $q$, $q$ odd and small, the size of the largest 2-intersecting family has been determined through computer search.
These values can be found in Table \ref{TableLaguerre2Intersecting}.
This table suggests that the true size is probably more or less linear in $q$, not quadratic.

\begin{table}[htbp]
    \centering
    \begin{tabular}{c|c c c c c c}
        $q$ & 3 & 5 & 7 & 9 & 11 & 13 \\ \hline
        size & 4 & 7 & 10 & 13 & 19 & 19 \\
        $\frac{q^2 - 1}2$ & 4 & 12 & 24 & 40 & 60 & 84
    \end{tabular}
    \caption{This table gives the size of the largest 2-intersecting families in an ovoidal Laguerre plane of order $q$, for $q \leq 13$ odd.
    We compare it with the bound of Theorem \ref{ThmLaguerre2intersecting}.}
    \label{TableLaguerre2Intersecting}
\end{table}

\subsection{The association scheme for ovoidal Laguerre planes of even order}
 \label{SubsectionLaguerreEvenAssoc}

Throughout the rest of this section, we will assume that $q > 2$ is even.
Let $L$ denote an ovoidal Laguerre plane of order $q$.
Let $\mq$ denote the oval cone in $\pg(3,q)$ from which $L$ is constructed.
Let $R$ be the vertex of this cone.
Take an oval plane $\pi$.
Since $\pi$ intersects $\mq$ in an oval and $q$ is even, this oval has a nucleus $N$.
We call $N$ the nucleus of $\pi$.
Take another oval plane $\rho$.
Then $\Psi_{\pi,\rho}$ as defined in Remark \ref{RmkIsomorphicOvalPlanes} is an collineation from $\pi$ to $\rho$, which maps $\pi \cap \mq$ to $\rho \cap \mq$.
Therefore, it must map the nucleus of $\pi$ to the nucleus of $\rho$.
This means that the nucleus of $\rho$ lies on the line $\nu = \vspan{N,R}$.
Thus, for every oval plane $\rho$, the nucleus of $\rho$ is $\rho \cap \nu$.

We make a slight modification to the incidence structure $L$ of the Laguerre plane.
Let $\mq_+$ denote the set $(\mq \cup \nu) \setminus \set R$.
Consider the incidence structure $L_+$ with $\mq_+$ as point set, and the oval planes as blocks.
(Technically, a block in this incidence structure is the intersection of $\mq_+$ and an oval plane.)
In the original Laguerre plane $L$, two oval planes $\pi$ and $\rho$ intersect in a unique point of $\mq$ if and only if the line $\pi \cap \rho$ is a tangent, and thus contains the nuclei of $\pi$ and $\rho$.
This happens if and only if $\pi$ and $\rho$ have the same nucleus.
Therefore, no two circles in $L_+$ intersect in exactly one point.
An intersecting family of circles in $L$ is the same thing as a 2-intersecting family in $L_+$.

We prove that the relations $R_0, R_2, R_3$ define an association scheme on the circles of $L_+$.
We do this be using a link between two-weight codes and strongly regular graphs, first explored by Delsarte \cite{delsarte71, delsarte72}.
We also refer to the survey of Calderbank and Kantor \cite{calderbankkantor}.

\begin{res}[{\cite[Theorem 3.2, Corollary 3.5]{calderbankkantor}}]
 \label{Res2WeightSRG}
Let $\mathcal A$ be a set in $\pg(k-1,q)$ that spans the entire space.
Let $\Omega$ denote the set $\sett{v \in \FF_q^k}{\vspan{v} \in \mathcal A}$.
Let $G(\Omega)$ denote the graph with as vertices the vectors of $\FF_q^k$, where two vertices $v$ and $w$ are adjacent if and only if $v-w \in \Omega$.
Then the following are equivalent:
\begin{enumerate}
 \item $|\mathcal A| = n$ and $\mathcal A$ intersects every hyperplane of $\pg(k-1,q)$ in $n-w_1$ or $n-w_2$ points.
 \item The adjacency matrix of $G(\Omega)$ has eigenvalues $n(q-1)$ with $\vspan{\one}$ as eigenspace, $n(q-1) - q w_1$, and $n(q-1) - q w_2$.
\end{enumerate}
\end{res}

We can choose coordinates of $\pg(3,q)$ in such a way that $R = (0,0,0,1)$.
The oval planes are exactly the planes missing $R$, hence the planes $\pi_{a,b,c}$ with equation $X_3 = a X_0 + b X_1 + c X_2$.
Now let $\mathcal O$ denote the intersection of $\mq_+$ and $\pi_{0,0,0}$.
Then $\mathcal O$ is a so-called hyperoval, i.e.\ a set of $q+2$ points in a projective plane of order $q$, such that each line intersects $\mathcal O$ in 0 or 2 points.
There are $\binom{q+2} 2$ secant lines and $\binom q 2$ skew lines.
Define the set
\[
 \Omega = \sett{(a,b,c) \in \FF_q^3}{ ( \, \forall (x,y,z,0) \in \mathcal O \,) (\, a x + b y + c z \neq 0 \,)}.
\]
Then a line in $\pi_{0,0,0}$ is skew to $\mathcal O$ if and only if this line has an equation of the form $X_3 = a X_0 + b X_1 + c X_2 = 0$ with $(a,b,c) \in \Omega$.
A point $P \in \mathcal O$ lies on no skew lines to $\mathcal O$, a point $P \in \pi_{0,0,0} \setminus \mathcal O$ lies on $\frac q 2$ skew lines to $\mathcal O$.
Dually, this means that if $\mathcal A$ denotes the set $\sett{\vspan v}{v \in \Omega}$ in $\pg(2,q)$, then every line in $\pg(2,q)$ intersects $\mathcal A$ in 0 or $\frac q 2$ points.
Thus, $\mathcal A$ and $\Omega$ are as in Result \ref{Res2WeightSRG}, with $n = \binom q 2$, $n - w_1 = \frac q 2$, and $n-w_2 = 0$.

Now let $G$ denote the graph with adjacency matrix $A_3$, i.e.\ the vertices of $G$ are the oval planes, and two vertices $\pi$ and $\rho$ are adjacent if and only if $\pi \cap \rho \cap \mq_+ = \emptyset$.

\begin{lm}
 \label{LmSRG}
The graph $G(\Omega)$ is isomorphic to $G$.
\end{lm}

\begin{proof}
We map the vertex $\pi_{a,b,c}$ of $G$ to the vertex $(a,b,c)$ of $G(\Omega)$.
To prove that this is a graph isomorphism, we need to show that $\pi_{a,b,c} \cap \pi_{\alpha,\beta,\gamma} \cap \mq_+ = \emptyset$ if and only if $(a,b,c) - (\alpha,\beta,\gamma) \in \Omega$.

So take two distinct oval planes $\pi_{a,b,c}$ and $\pi_{\alpha,\beta,\gamma}$.
They intersect in the line $l$ with equation $X_3 = a X_0 + b X_1 + c X_2$ and $(a-\alpha) X_0 + (b - \beta) X_1 + (c - \gamma) X_2 = 0$.
Note that $l \cap \mq_+ = \emptyset$ if and only if $\vspan{l, R} \cap \mq_+ = \emptyset$.
The plane $\vspan{l,R}$ has equation $(a-\alpha) X_0 + (b - \beta) X_1 + (c - \gamma) X_2 = 0$.
The set $\mq_+$ consists of all the points $(x,y,z,u)$ with $(x,y,z,0) \in \mathcal O$ and $u \in \FF_q$.
Thus, no point of $\mq_+$ lies in $\vspan{l,R}$ if and only if there are no solutions to $(a-\alpha) x + (b - \beta) y + (c-\gamma) z = 0$ with $(x,y,z,0) \in \mathcal O$.
This is equivalent to $(a,b,c) - (\alpha,\beta,\gamma) \in \Omega$.
\end{proof}

Hence, we can use Result \ref{Res2WeightSRG} to calculate the spectrum of $A_3$.
In the notation of Result \ref{Res2WeightSRG}, $n = w_2= \binom q 2$, $w_1 = (q-2)\frac q 2$, which means that the eigenvalues of $A_3$ are $n(q-1) = (q-1)^2 \frac q 2$, $n(q-1) - q w_1 = \frac q 2$, and $n(q-1) - q w_2 = - (q-1) \frac q 2$.
Regular graphs with only three distinct eigenvalues can be found in the literature under the name \emph{strongly regular graphs}, see e.g.\ \cite[\S 5]{godsilmeagher}, \cite[\S 1.3]{bcn}.
There is a simple and well-known way to calculate the multiplicities of the eigenvalues.
Let $m_1$ and $m_2$ denote the multiplicities of $\frac q 2$ and $-(q+1)\frac q 2$ respectively.
The sum of the multiplicities of the eigenvalues of $A_3$ equals the size of $A_3$, hence $1 + m_1 + m_2 = q^3$.
The sum of the eigenvalues of $A_3$ counted with multiplicity equals the trace of $A_3$, hence $1 \cdot (q-1)^2 \frac q 2 + m_1 \frac q 2 - m_2 (q-1) \frac q 2 = 0$.
This implies that $m_1 = (q+1)(q-1)^2$ and $m_2 = (q-1)(q+2)$.
Then $A_2 = J - I - A_3$, and it is easy to check that $I, A_2, A_3$ constitute a 2-class association scheme.
Let $V_1$ and $V_2$ denote the eigenspace of $A_3$ for eigenvalue $\frac q 2$ and $-(q+1)\frac q 2$ respectively.
Then $A_2$ has eigenvalue $\binom{q+2}2 (q-1)$ on $\vspan{\one}$, $-\frac{q+2}2$ on $V_1$, and $(q^2+1) \frac{q-2}2$ on $V_2$.

\subsection{Intersecting families in ovoidal Laguerre planes of even order}

Let $G$ still denote the graph with adjacency matrix $A_3$.
Then 2-intersecting families in $L_+$ are equivalent to cocliques in $G$.
We apply Hoffman's ratio bound.
\[
 \alpha(G) \leq \frac{q^3}{ 1 + \frac{ {(q-1)^2} \frac q2 }{ (q-1)\frac q 2 }} = q^2.
\]
A coclique attaining this bound has its characteristic vector in $\vspan{\one, V_2}$.
Consider the incidence matrix $W$ of $L_+$.
Then $W W^t$ is a symmetric matrix, indexed by the oval planes, and $W W^t (\pi,\rho)$ equals the number of points in $\pi \cap \rho \cap \mq_+$.
Therefore,
\[
 W W^t = (q+2) I + 2 A_2.
\]
Hence, $\ker(W W^t) = V_1$, and $\col(W) = \col(W W^t) = \ker(W W^t)^\perp = \vspan{\one,V_2}$.

Thus, if $K$ is 2-intersecting family in $L_+$ of size $q^2$, then its characteristic vector equals $W \kappa$ for some vector $\kappa$.
As was done previously, take an oval plane $\pi \in K$.
Let $W_\pi$ denote the restriction of $W$ to the rows corresponding to oval planes in relation $R_3$ with $\pi$ and columns corresponding to the points of $\mq_+ \setminus \pi$.
Let $\kappa_\pi$ denote the restriction of $\kappa$ to the points of $\mq_+ \setminus \pi$.
Then $W_\pi \kappa_\pi = 0$.
We determine $\ker(W_\pi) = \ker(W_\pi^t W_\pi)$.
The $(P,Q)$-entry of $W_\pi^t W_\pi$ equals the number of oval planes through $P$ and $Q$, intersecting $\pi$ in a skew line.
We extend the parallel relation of $L$ by taking $\nu \setminus \set R$ as parallel class.
In other words, any two points $P$ and $Q$ of $\mq_+$ are parallel if and only if $P$, $Q$, and $R$ are collinear.

First suppose that $P=Q$.
Each of the $(q-1)\frac q 2$ skew lines in $\pi$ lies in a unique plane through $P$.
Every plane through $P$ and $R$ contains the point parallel to $P$ in $\pi$, and hence doesn't intersect $\pi$ in a skew line.
Therefore, there are $(q-1) \frac q 2$ oval planes through $P$ intersecting $\pi$ in a skew line.

Secondly suppose that $P\neq Q$, but $P$ and $Q$ are parallel.
Then no plane through $P$ and $Q$ is an oval plane.

Thirdly, suppose that $P$ and $Q$ are not parallel.
Consider the line $m = \vspan{P,Q}$.
Since $P$ and $Q$ are not parallel, $R \notin m$.
Then $m$ can contain at most two points of $\mq_+$, thus $m \cap \mq_+ = \set{P,Q}$.
In particular, $m \cap \pi \notin \mq_+$.
Thus, $m \cap \pi$ is a point outside of $\mq_+$, and therefore lies on $\frac q 2$ skew lines in $\pi$, each giving us a unique plane through $m$.
All these planes are oval planes, since there is only one plane through $m$ and $R$, and this plane contains the points of $\pi \cap \mq_+$ parallel to $P$ and $Q$.
We may conclude that
\[
 W_\pi^t W_\pi (P,Q) = \begin{cases}
  (q-1)\frac q 2 & \text{if } P=Q, \\
  0 & \text{if $P \neq Q$, and $P$ and $Q$ are parallel}, \\
  \frac q 2 & \text{if $P$ and $Q$ aren't parallel.}
 \end{cases}
\]
Suppose that the ordering of the columns of $W_\pi$ consists of $q+2$ consecutive blocks, each of $q-1$ parallel points.
Then
\[
 W_\pi^t W_\pi = \frac q 2 \Big( (q-1)  I_{(q+2)(q-1)} + (J_{q+2}-I_{q+2}) \otimes J_{q-1} \Big).
\]
Hence, the kernel of $W_\pi^t W_\pi$ equals of the eigenspace of $(J_{q+2}-I_{q+2}) \otimes J_{q-1}$ with eigenvalue $-(q-1)$.
The matrix $J_{q+2}-I_{q+2}$ has eigenvalue $q+1$ on $\vspan{\one_{q+1}}$, and eigenvalue $-1$ on $\vspan{\one_{q+2}}^\perp$.
The matrix $J_{q-1}$ has eigenvalue $q-1$ on $\vspan{\one_{q-1}}$, and eigenvalue $0$ on $\vspan{\one_{q-1}}^\perp$.
As can be deduced from \S \ref{SubsectionKronecker}, the eigenvectors of $(J_{q+2}-I_{q+2}) \otimes J_{q-1}$ with eigenvalue $-(q-1)$ are of the from $v \otimes \one_{q-1}$, with $v \in \vspan{\one_{q+2}}^\perp$.
Thus, $\kappa_\pi = v_\pi \otimes \one_{q-1}$ for some $v_\pi \in \vspan{\one_{q+1}}^\perp$.
This means the following.
We can interpret $v_\pi$ as a function on the parallel classes.
Recall that we denote the parallel class containing a point $P$ as $\overline P$.
For each point $P \in \mq_+ \setminus \pi$, $\kappa(P) = v_\pi(\overline P)$, and $\sum_{P \in \pi \cap \mq_+} v_\pi (\overline P) = 0$.
Such a vector $v_\pi$ exists for each $\pi \in K$.

The end of the proof is analogous to the end of the proof for ovoidal Laguerre planes of odd order.
Take two planes $\pi$ and $\rho$ in $K$.
Take a point $P \in \pi \cap \mq$, and let $Q$ be the point parallel with $P$ in $\rho$.
Since we assumed that $q>2$, there exists a point $T \in \overline P \setminus \set{P,Q}$.
Then $v_\pi(\overline P) = \kappa(T) = v_\rho(\overline P)$.
Thus, $v_\pi$ is equal for all $\pi \in K$.
We drop the subscript and denote this vector as $v$.
If for each point $P \in \mq_+$ there exists a plane $\pi \in K$ with $P \notin \pi$, then $\kappa(P) = v(\overline P)$ for all points $P$.
Take the row $w$ of $W$ corresponding to the oval plane $\rho$.
Then $w \kappa = \sum_{P \in \rho} \kappa(P) = \sum_{P \in \rho} v(\overline P) = 0$ since $v \in \vspan{\one}^\perp$.
But then $W \kappa = 0$, which means that $K = \emptyset$.
This is a contradiction, thus there must exist a point $P$ that lies in all planes of $K$.
Since $|K| = q^2$ equals the number of oval planes through a point, $K$ consists exactly of all oval planes through $P$.

Therefore, the 2-intersecting families of size $q^2$ in $L_+$ are exactly the families consisting of all oval planes through a fixed point.
This means that there are two types of intersecting families in the Laguerre plane $L$.

\begin{thm}
 Let $\mathcal F$ be an intersecting family in an ovoidal Laguerre plane $L = (\mathcal P, \mathcal B)$ of even order $q > 2$.
 Then $|\mathcal F| \leq q^2$.
 Equality holds if and only if $\mathcal F$ consists of all circles through a fixed point of $\mathcal P$, or all circles with a fixed nucleus.
\end{thm}

Since two circles in $L$ intersect in a unique point of $\mq$ if and only if they share a nucleus, this theorem is clearly equivalent to Theorem \ref{ThmMainLaguerreEven}.

\bigskip 

What can we say for $q=2$?
The only ovals in $\pg(2,2)$ are the conics, so we can use the polynomial description for the Laguerre plane of order 2.
Two polynomials $f$ and $g$ are in relation $R_3$ if and only if $f-g$ has no roots in $\FF_2 \cup \set \infty$.
The only polynomial of degree at most 2 with no roots over $\FF_2 \cup \set \infty$ is $X^2 + X + 1$.
Thus, an intersecting family is the same thing as a set of polynomials $K$ such that for any $f \in K$, $X^2 + X + 1 - f \notin K$.
Therefore, to make an intersecting family, we need to choose at most 1 polynomial from each of the 4 pairs of the form $\set{f,X^2+X+1-f}$.
This way we see that the maximum size is $4 = q^2$.
However, there are $2^4 = 16$ options to create such a family.
On the other hand, there are only 6 points in $(\FF_2 \cup \set \infty) \times \FF_2$, and 2 nuclei in the geometrical description of the Laguerre plane.
Hence, the ``classical'' intersecting families only account for 8 of the 16 intersecting families of size 4.
This gives us a weak, but definitely not a strong EKR property.

\subsection{2-intersecting families in ovoidal Laguerre planes of even order}
 \label{SubsectionLaguerre2intersectingEven}
2-intersecting families in ovoidal Laguerre planes of odd order satisfy the weak EKR-property, which can be proven very easily.

\begin{thm}
 Let $\mathcal F$ be a 2-intersecting family in an ovoidal Laguerre plane of even order $q$.
 Then $|\mathcal F| \leq q$, and this bound is tight.
\end{thm}

\begin{proof}
If two oval planes share a nucleus, they are intersecting.
Since there are only $q$ nuclei, a 2-intersecting family has at most $q$ members.
Furthermore, given two non-singular points $P$ and $Q$ of an oval cone, there are $q$ oval planes through the line $\vspan{P,Q}$, which yields a 2-intersecting family of size $q$.
\end{proof}

\section{Polynomials of bounded degree} 
 \label{SectionPolynomials}
Consider the vector space $U_{q,k} = \sett{f \in \FF_q[X]}{\deg f \leq k}$ of polynomials over $\FF_q$ of degree at most $k$.
We will denote it by $U$ if $q$ and $k$ are clear from context.
In this section, we investigate whether the EKR properties are satisfied in the incidence structure with $\FF_q \times \FF_q$ as points, and the graphs of the polynomials in $U$ as blocks.

Say that two polynomials $f$ and $g$ in $U$ are $t$-\emph{intersecting} if there exist at least $t$ values $x \in \FF_q$ such that $f(x) = g(x)$, or equivalently if $f-g$ has at least $t$ distinct roots.
Call a family $\mathcal F$ of polynomials $t$-intersecting if any two polynomials in $\mathcal F$ are $t$-intersecting.
Call $\mathcal F$ non-$t$-intersecting if no two distinct elements of $\mathcal F$ are $t$-intersecting.
We are interested in the largest set of pairwise $t$-intersecting polynomials.

\begin{lm}
 \label{LmTIntPol}
Assume that $t \leq k < q$.
The largest $t$-intersecting family in $U_{q,k}$ has size $q^{k+1-t}$.
The largest non-$t$-intersecting family in $U_{q,k}$ has size $q^t$.
\end{lm}

\begin{proof}
Consider the graph $G$ with $U_{q,k}$ as vertices, and where $f$ and $g$ are adjacent if and only if $f$ and $g$ are $t$-intersecting.
Cliques and cocliques in this graph are exactly the $t$-intersecting and non-$t$-intersecting families of $U$ respectively.
Thus, we need to prove that $\omega(G) = q^{k-t+1}$ and that $\alpha(G) = q^t$.
Choose $t$ distinct values $x_1, \dots, x_t$ in $\FF_q$, and $t$ (not necessarily distinct) numbers $y_1, \dots, y_t$ in $\FF_q$.
Then all polynomials $f$ such that $(\, \forall i \in [1,t]\, )(\, f(x_i)=y_i \, )$ are a set of $q^{k+1-t}$ pairwise $t$-intersecting polynomials, hence $\omega(G) \geq q^{k-t+1}$.
On the other hand, $U_{q,t-1}$ is a non-$t$-intersecting family of size $q^t$, hence $\alpha(G) \geq q^t$.
We can easily prove that $G$ is vertex transitive, since it is a Cayley graph:
for each $f \in U$, the map $U \to U: g \mapsto g+f$ is an automorphism of $G$, and these automorphisms work sharply transitive on the vertices of $G$.
Furthermore, $G$ has $|U| = q^{k+1}$ vertices, hence the clique-coclique bound implies that $\omega(G) \alpha(G)$ is at most $q^{k+1}$, which finishes the proof.
\end{proof}

This proves the weak EKR-property for $t$-intersecting families of polynomials.
For 1-intersecting families in $U_{q,k}$, $2 \leq k < q$, the strong EKR-property holds.
To state the theorem, let $\mathcal F_{x,y}$ denote the set $\sett{f \in U_{q,k}}{f(x) = y}$, $x,y \in \FF_q$.

\begin{thm}
 \label{ThmPolynomialsIntersectingEKR}
Assume that $2 \leq k < q$.
Then the intersecting families of size $q^k$ in $U_{q,k}$ are exactly the families $\mathcal F_{x,y}$, for $x,y \in \FF_q$.
\end{thm}

\begin{proof}
The theorem will be proven through induction on $k$.
First assume that $k=2$.
Let $\mathcal F$ be an intersecting family in $U_{q,2}$ of size $q^2$.
If we extend each function $f: \FF_q \to \FF_q: X \mapsto aX^2 + bX + c$ to a function on $\FF_q \cup \set \infty$ by defining $f(\infty) = a$, then the polynomials of $U_{q,2}$ represent the circles of an ovoidal Laguerre plane (see Remark \ref{RmkOvoidalLaguerrePolynomials}.)
Note that in this interpretation $\mathcal F$ is a still an intersecting family.
Suppose that $q$ is odd.
By Theorem \ref{ThmMainLaguerreOdd}, the strong EKR property holds and $\mathcal F = \mathcal F_{x,y}$, for some $x \in \FF_q \cup \set \infty$, and $y \in \FF_q$.
If $x = \infty$, then $\mathcal F$ contains $yX^2$ and $yX^2 +1$, which are not intersecting in our original interpretation.
Thus, $x \neq \infty$.
Now suppose that $q$ is even.
The assumptions from the theorem imply that $q>2$.
By Theorem \ref{ThmMainLaguerreEven}, $\mathcal F = \mathcal F_{x,y}$, for some $x \in \FF_q \cup \set \infty$, and $y \in \FF_q$, or $\mathcal F$ consists of a polynomial $f$ and all polynomials which are exactly 1-intersecting with $f$ in the Laguerre interpretation.
Two polynomials belong to a family of the latter type if and only if their linear coefficient is equal.
This means that such a family contains polynomials $f$ and $f+1$, which aren't intersecting in our original interpretation.
Thus $\mathcal F$ is of the former type.
We can exclude $x = \infty$ as we did for $q$ odd.
Thus, the case $k=2$ is proven.

Assume that $U_{q,k}$ satisfies the strong EKR property, $k \leq q-2$.
This assumption is the induction hypothesis.
Take an intersecting family $K$ in $U = U_{q,k+1}$ of size $q^{k+1}$.
Divide $K$ into $q$ classes $K_\alpha$, $\alpha \in \FF_q$, where $f \in K_\alpha$ if and only if the leading coefficient of $f$ is $\alpha$.
Note that $K_\alpha$ is isomorphic to $U_{q,k}$ with respect to the property of being $t$-intersecting (for any value of $t$).
Therefore, since $K_\alpha$ is an intersecting family, $|K_\alpha| \leq q^k$ by Lemma \ref{LmTIntPol}.
Since $|K| = q^{k+1}$, this implies that $|K_\alpha| = q^k$ for each $\alpha$.
Thus we know from the induction hypothesis that for every $\alpha$ there exist $x_\alpha$ and $y_\alpha$ such that $f(x_\alpha)=y_\alpha$ for all $f$ in $K_\alpha$.

\underline{Case 1: $x_0 = y_0 = 0$.}

We want to show that $x_\alpha = y_\alpha=0$ for all other $\alpha$'s.
In other words, given $f \in K_\alpha$, we must show that the constant term of $f$ is zero.
Assume by contradiction that there is a polynomial $f \in K_\alpha$ with
$f(X) = \alpha X^{k+1} + \sum_{i=1}^k \beta_i X^k + \gamma$, and $\gamma \neq 0$.
We prove that there exists a polynomial $g \in K_0$ not intersecting $f$.
Note that the set $\sett{f-g}{g\in K_0}$ are all polynomials with leading coefficient $\alpha$ and constant term $\gamma$.
Thus, fixing $\alpha$ and $\gamma$, it suffices to show that we can choose the values $\beta_i$ such that $f$ has no roots over $\FF_q$.
It follows from the theory of Vandermonde matrices that if we choose distinct non-zero values $x_1,\dots,x_j$ in $\FF_q^*$, $j\leq k$, then the matrix
\[
 V = \begin{pmatrix}
  1 & 0 & \dots & 0 & 0\\
  0 & 0 & \dots & 0 & 1\\
  x_1^{k+1} & x_1^k & \dots & x_1 & 1 \\
  \vdots &           &     &      & \vdots \\
  x_j^{k+1} & x_j^k & \dots & x_j & 1
 \end{pmatrix}
\]
has full rank.
This means that given $\alpha$ and $\gamma$ and predetermined values of $f$ in $x_1, \dots, x_j$, there are $q^{k-j}$ ways to choose the $\beta_i$ coefficients.
Given $k+1$ values $x_1, \dots, x_{k+1}$, there is a unique polynomial with leading coefficient $\alpha$ and $x_1, \dots, x_{k+1}$ as roots, namely $\alpha(X-x_1) \dots (X-x_{k+1})$.
Its constant coefficient is $\alpha (-1)^{k+1} x_1 \cdots x_{k+1}$.

Now define for each $x \in \FF_q^*$ the set $N_x$ of polynomials $f$ with leading coefficient $\alpha$ and constant term $\gamma$, for which $f(x)=0$.
By previous considerations, it follows that given distinct values $x_1, \dots, x_j \in \FF_q^*$,
\[
 \left| \bigcap_{1 \leq i \leq j} N_{x_i} \right|
 = \begin{cases}
  q^{k-j} & \text{if } j \leq k, \\
  1 & \text{if } j=k+1 \text{ and } x_1 \cdot \ldots \cdot x_{k+1} = (-1)^{k+1} \frac \gamma \alpha, \\
  0 & \text{otherwise}.
 \end{cases}
\]
Let $\binom{\FF_q^*}j$ denote all subsets of $\FF_q^*$ of size $j$.
By the inclusion-exclusion principle,
 \begin{align}
  \label{EqNx}
  \left| \bigcup_{x \in \FF_q^*} N_x \right|
  = \sum_{j=1}^{k+1} (-1)^{j-1} \sum_{A \in \binom{\FF_q^*}{j}} \left| \bigcap_{x\in A} N_x \right|
  = \sum_{j=1}^{k} (-1)^{j-1} \binom{q-1}{j} q^{k-j} + (-1)^k s,
 \end{align}
 where $s$ is the number of $(k+1)$-sets of $\FF_q^*$ for which the product of the elements equals $(-1)^{k+1} \alpha/\gamma$.
 Note that $0 \leq s \leq \frac 1 {k+1} \binom{q-1}{k}$, since each $k$-set can be extended in at most one way to a $(k+1)$-set with prescribed product, and every $(k+1)$-set has $k+1$ subsets of size $k$.
 Also note that
 \[
 \frac{\binom{q-1} j q^{k-j}}{\binom{q-1} {j+1} q^{k-(j+1)}}
 = \frac{ \frac{q}{j! (q-1-j)!}  }{ \frac{1}{(j+1)! (q-2-j)!}  }
 = \frac{q (j+1)}{q-1-j} \geq j+1.
 \]
 Hence, the right hand side of (\ref{EqNx}) is an alternating series where the terms decrease in absolute value.
 Thus, we can bound it from above by the first term $(q-1)q^{k-1}<q^k$.
 Therefore, there exists a polynomial $f$ with leading coefficient $\alpha$ and constant term $\gamma$, which is not a member of $\bigcup_{x \in \FF_q^*} N_x$, meaning $f$ has no roots.
 This completes the contradiction.

\underline{Case 2: $x_0$ and $y_0$ are not both zero.}

Consider the bijection $\sigma: U \to U: f \mapsto f(X+x_0) - y_0$.
Note that $\sigma$ maps $K$ to an intersecting set where $x_0 = y_0 = 0$.
Therefore, as was proven in Case 1, $x_\alpha = y_\alpha = 0$ for all $(\sigma(K))_\alpha$.
Reversing $\sigma$, this means that $x_\alpha=x_0$ and $y_\alpha=y_0$ for all $\alpha$, proving the strong EKR property.
\end{proof}

Hence, we have proven Theorem \ref{ThmMainPolynomials}.

\section{Minkowski planes}
 \label{SectionMinkowski}

\subsection{Intersecting families in ovoidal Minkowski planes}

An ovoidal Minkowski plane arises from the points and oval planes of a quadratic set in $\pg(3,q)$, with no singular points, containing a line.
By Buekenhout's result \cite{buekenhout69}, the only such sets are the hyperbolic quadrics (see \S \ref{SubsectionQuadraticSets}).
The corresponding Minkowski plane can be represented as the graphs of the elements of $\text{PGL}(2,q)$ \cite[\S 5]{hartmann}.
In other words, as the incidence structure $(\mathcal P, \mathcal B)$ with $\mathcal P = \pg(1,q) \times \pg(1,q)$, and the elements of $\mathcal B$ are the graphs $\sett{(P,f(P))}{P \in \pg(1,q)}$ of $f \in \text{PGL}(2,q)$.
As a consequence, the main theorem of \cite{meagherspiga} is equivalent to the following theorem.

\begin{thm}
 \label{ThmMinkowski}
Let $M = (\mathcal P, \mathcal B)$ be an ovoidal Minkowski plane of order $q$.
Let $\mathcal F$ be an intersecting family in $M$.
Then $|\mathcal F| \leq q(q-1)$.
Equality occurs if and only if $\mathcal F$ consists of all circles through a point of $\mathcal P$.
\end{thm}

For ovoidal Minkowski planes of even order, the relations of Definition \ref{DfRelAssoc} constitute an association scheme.
This does not hold for odd values of $q$, because the graph corresponding to relation $R_3$ has 5 distinct eigenvalues (see \cite[Table 3]{meagherspiga}).
For the sake of completion, we mention that the matrices of eigenvalues and dual eigenvalues of the association scheme for $q$ even are
\begin{align*}
 P = \begin{pmatrix}
 1 & q^2-1 & q(q+1)\frac{q-2}2 & (q-1)\frac{q^2}2 \\
 1 & q-1 & -q & 0 \\
 1 & -(q+1) & 0 & q \\
 1 & 0 & \frac{q^2-q-2}2 & - (q-1)\frac q 2
 \end{pmatrix}, &&
 Q = \begin{pmatrix}
 1 & (q+1)^2\frac{q-2}2 & (q-1)^2\frac q 2 & q^2 \\
 1 & (q+1)\frac{q-2}2 & -(q-1) \frac q 2 & 0 \\
 1 & -(q+1) & 0 & q \\
 1 & 0 & q-1 & -q
 \end{pmatrix}.
\end{align*}
The proof that these are the correct matrices works similar to the previous ones, by calculating the intersection numbers and the intersection matrices.
Since we will not really use the association scheme here, we omit the proof from this paper.

\subsection{2-intersecting families in Minkowski planes}

Using spectral techniques for bounding the size of 2-intersecting families an ovoidal Minkowski plane of even order $q$, yields $1 + (q+1)\frac{q-2}2$.
Again, this can also be proven by a Blokhuis-Bruen type argument, as was the case for Theorem \ref{ThmLaguerre2intersecting}.

\bigskip

Take a 2-intersecting family $K$ is a Minkowski plane of order $q$.
Take a point $P$, and look in the residue at $P$.
This residue is an affine plane of order $q$ in which two parallel classes correspond to the parallel classes of the Minkowski plane.
Therefore, there are at most $q-1$ circles of $K$ through $P$.
Double count the set $Z$ as we did for Theorem \ref{ThmLaguerre2intersecting}.
This yields the bound $|K| \leq 1 + (q+1)\frac{q-2}2$.

We can improve this bound if $q>2$.
If it were tight, then every point is covered by $q-1$ or 0 circles of $K$.
This implies that the circles of $K$ cover $(q+1)\frac q 2$ points (perform a double count).
Assume that there are $q-1$ circles of $K$ through $P$.
How many points do they contain?
Interpret them as lines in the affine residue through $P$.
These circles then become lines pairwise intersecting in a point.
Denote these lines as $l_1, \dots, l_{q-1}$.
Then $l_1$ covers $q$ points, adding $l_2$ gives us $q-1$ extra points, adding $l_3$ gives us at least $q-2$ extra points etc.
In total, they cover at least $q + (q-1) + \dots + 2 = (q-1)\frac{q+2}2$ points.
Note that $(q-1)\frac{q+2}2 = (q+1)\frac q 2 - 1$.
Thus, all points covered by the circles in $K$ distinct from $P$, must lie in the residue at $P$.
This means that no point parallel to $P$ is covered.
Therefore, if every point lies on 0 or $q-1$ circles of $K$, then the points covered by the circles of $K$ are pairwise not parallel, which implies that there are at most $q+1$ such points.
This gives a contradiction since $q+1 < (q+1)\frac q 2$ if $q>2$.

If $q=2$, the bound states $|K| \leq 1$, which is evidently tight.

\begin{thm}
 \label{ThmMinkowski2intersecting}
If $\mathcal F$ is a 2-intersecting family in a Minkowski plane of order $q>2$, then $|\mathcal F| \leq (q+1)\frac{q-2}2$.
\end{thm}

A computer search found the size of the largest 2-intersecting families in the ovoidal Minkowski planes of order $q$ for small values of $q$, see Table \ref{TableMinkowski2Intersecting}.
The correct size seems to be more or less linear in $q$, in contrast to the quadratic bound from Theorem \ref{ThmMinkowski2intersecting}.

\begin{table}[htbp]
    \centering
    \begin{tabular}{c|c c c c c c c c c c c}
        $q$ & 2 & 3 & 4 & 5 & 7 & 8 & 9 & 11 & 13 & 16 & 17  \\ \hline
        size & 1 & 2 & 4 & 5 & 8 & 10 & 12 & 17 & 17 & 28 & 23 \\
        $(q+1)\frac{q-2}2$ &  & 2 & 5 & 9 & 20 & 27 & 35 & 54 & 77 & 119 & 135
    \end{tabular}
    \caption{This table gives the size of the largest 2-intersecting family in the ovoidal Minkowski plane of order $q \leq 17$.
    We compare it with the bound from Theorem \ref{ThmMinkowski2intersecting}, which holds for $q>2$.}
    \label{TableMinkowski2Intersecting}
\end{table}

\section{Concluding remarks}
 \label{SectionGeMoogtNaarHuisGaan}

This article can serve as an illustration of the potency of the linear algebraic approach to solving EKR-type problems.
Some open problems remain.
\begin{itemize}
\item Do intersecting families satisfy the strong EKR property in ovoidal Möbius planes of odd order?
\item What can we say about non-ovoidal circle geometries?
\item What are the largest 2-intersecting families in circle geometries?
\item Can we classify the largest $t$-intersecting families of polynomials in $U_{q,k}$?
\end{itemize}

\bigskip

It is remarkable that the algebraic approach is very fruitful to study 1-intersecting families in ovoidal circle geometries, but is easily surpassed by elementary combinatorial arguments to study 2-intersecting families.

\bigskip

We also mention that there is a natural correspondence between $U_{q,k}$ and Reed-Solomon codes over $\FF_q$ of dimension $k+1$ and length $q$.
Therefore, the results in this paper can be translated to results on Reed-Solomon codes.
The classical Laguerre planes corresponds to the extended Reed-Solomon codes of dimension $3$.
This raises the question whether we can say something in general about the EKR properties for MDS codes.

\bigskip

{\bf Acknowledgements.} The author would like to thank his supervisor Jan De Beule, whose fluency in GAP is a great convenience.
The author would also like to thank Sam Mattheus and Karen Meagher for helpful discussions.

\bibliographystyle{alpha}
\bibliography{ErdosKoRadoTheoremsForOvoidalCircleGeometries.bib}

\begin{thebibliography}{MST16}

\bibitem[AK97]{ahlswedekhachatrian}
R.~Ahlswede and L.~H. Khachatrian.
\newblock The complete intersection theorem for systems of finite sets.
\newblock {\em European J. Combin.}, 18(2):125--136, 1997.

\bibitem[Bar55]{barlotti55}
A.~Barlotti.
\newblock Un'estensione del teorema di {S}egre-{K}ustaanheimo.
\newblock {\em Boll. Un. Mat. Ital. (3)}, 10:498--506, 1955.

\bibitem[Bar65]{barlotti65}
A.~Barlotti.
\newblock Some topics in finite geometrical structures.
\newblock Technical report, North Carolina State University. Dept. of
  Statistics, 1965.

\bibitem[BB89]{blokhuisbruen}
A.~Blokhuis and A.~A. Bruen.
\newblock The minimal number of lines intersected by a set of {$q+2$} points,
  blocking sets, and intersecting circles.
\newblock {\em J. Combin. Theory Ser. A}, 50(2):308--315, 1989.

\bibitem[BCN89]{bcn}
A.~E. Brouwer, A.~M. Cohen, and A.~Neumaier.
\newblock {\em Distance-Regular Graphs}.
\newblock Springer, 1989.

\bibitem[Ben73]{benz73}
W.~Benz.
\newblock {\em Vorlesungen \"{u}ber {G}eometrie der {A}lgebren}.
\newblock Springer-Verlag, Berlin-New York, 1973.
\newblock Geometrien von M\"{o}bius, Laguerre-Lie, Minkowski in einheitlicher
  und grundlagengeometrischer Behandlung, Die Grundlehren der mathematischen
  Wissenschaften, Band 197.

\bibitem[Bue69]{buekenhout69}
F.~Buekenhout.
\newblock Ensembles quadratiques des espaces projectifs.
\newblock {\em Math. Z.}, 110:306--318, 1969.

\bibitem[Che70]{chen70}
Y.~Chen.
\newblock Der {S}atz von {M}iquel in der {M}\"{o}biusebene.
\newblock {\em Math. Ann.}, 186:81--100, 1970.

\bibitem[Che74]{chen74}
Y.~Chen.
\newblock A characterization of some geometries of chains.
\newblock {\em Canadian J. Math.}, 26:257--272, 1974.

\bibitem[CK86]{calderbankkantor}
R.~Calderbank and W.~M. Kantor.
\newblock The geometry of two-weight codes.
\newblock {\em Bull. London Math. Soc.}, 18(2):97--122, 1986.

\bibitem[CK03]{cameronku}
P.~J. Cameron and C.~Y. Ku.
\newblock Intersecting families of permutations.
\newblock {\em European Journal of Combinatorics}, 24(7):881--890, 2003.

\bibitem[DB14]{maartenthesis}
M.~De~Boeck.
\newblock {\em Intersection problems in finite geometries}.
\newblock PhD thesis, Ghent University, 2014.

\bibitem[Del71]{delsarte71}
Ph. Delsarte.
\newblock {\em Two-weight linear codes and strongly regular graphs}.
\newblock MBLE. Laboratoire de Recherches, 1971.

\bibitem[Del72]{delsarte72}
Ph. Delsarte.
\newblock Weights of linear codes and strongly regular normed spaces.
\newblock {\em Discrete Math.}, 3:47--64, 1972.

\bibitem[Del73]{delsartethesis}
Ph. Delsarte.
\newblock An algebraic approach to the association schemes of coding theory.
\newblock {\em Philips Res. Rep. Suppl.}, 10:vi+97, 1973.

\bibitem[Dem64]{dembowski64}
P.~Dembowski.
\newblock M\"{o}biusebenen gerader {O}rdnung.
\newblock {\em Math. Ann.}, 157:179--205, 1964.

\bibitem[Dem68]{dembowksi68}
P.~Dembowski.
\newblock {\em Finite geometries}.
\newblock Ergebnisse der Mathematik und ihrer Grenzgebiete, Band 44.
  Springer-Verlag, Berlin-New York, 1968.

\bibitem[EKR61]{erdoskorado}
P.~Erd\H{o}s, C.~Ko, and R.~Rado.
\newblock Intersection theorems for systems of finite sets.
\newblock {\em Quart. J. Math. Oxford Ser. (2)}, 12:313--320, 1961.

\bibitem[FD77]{frankldeza}
P.~Frankl and M.~Deza.
\newblock On the maximum number of permutations with given maximal or minimal
  distance.
\newblock {\em Journal of Combinatorial Theory, Series A}, 22(3):352--360,
  1977.

\bibitem[GM16]{godsilmeagher}
C.~Godsil and K.~Meagher.
\newblock {\em Erd\H{o}s-{K}o-{R}ado theorems: algebraic approaches}, volume
  149 of {\em Cambridge Studies in Advanced Mathematics}.
\newblock Cambridge University Press, Cambridge, 2016.

\bibitem[God10]{godsil}
C.~Godsil.
\newblock Association schemes.
\newblock \url{https://www.math.uwaterloo.ca/~cgodsil/pdfs/assoc2.pdf}, 2010.

\bibitem[Hae21]{haemers2021}
W.~Haemers.
\newblock Hoffman's ratio bound.
\newblock \url{https://arxiv.org/abs/2102.05529}, 2021.

\bibitem[Har04]{hartmann}
E.~Hartmann.
\newblock Planar circle geometries.
\newblock
  \url{https://www2.mathematik.tu-darmstadt.de/~ehartmann/circlegeom.pdf},
  2004.

\bibitem[Hsi75]{hsieh75}
W.N. Hsieh.
\newblock Intersection theorems for systems of finite vector spaces.
\newblock {\em Discrete Mathematics}, 12(1):1--16, 1975.

\bibitem[Kah80]{kahn80}
J.~Kahn.
\newblock Locally projective-planar lattices which satisfy the bundle theorem.
\newblock {\em Math. Z.}, 175(3):219--247, 1980.

\bibitem[MS11]{meagherspiga}
K.~Meagher and P.~Spiga.
\newblock An {E}rd{\H{o}}s-{K}o-{R}ado theorem for the derangement graph of
  {${\rm PGL}(2,q)$} acting on the projective line.
\newblock {\em J. Combin. Theory Ser. A}, 118(2):532--544, 2011.

\bibitem[MST16]{meagherspigatiep}
K.~Meagher, P.~Spiga, and P.~H. Tiep.
\newblock An erd{\H{o}}s--ko--rado theorem for finite 2-transitive groups.
\newblock {\em European Journal of Combinatorics}, 55:100--118, 2016.

\bibitem[PSV11]{pepestormevanhove}
V.~Pepe, L.~Storme, and F.~Vanhove.
\newblock Theorems of erd{\H{o}}s--ko--rado type in polar spaces.
\newblock {\em Journal of Combinatorial Theory, Series A}, 118(4):1291--1312,
  2011.

\bibitem[Seg55]{segre55}
B.~Segre.
\newblock Ovals in a finite projective plane.
\newblock {\em Canadian J. Math.}, 7:414--416, 1955.

\bibitem[SS97]{willihans}
W.-H. Steeb and T.~K. Shi.
\newblock {\em Matrix calculus and Kronecker product with applications and C++
  programs}.
\newblock World Scientific, 1997.

\bibitem[Tit60]{tits60}
J.~Tits.
\newblock Les groupes simples de suzuki et de ree.
\newblock {\em S{\'e}minaire Bourbaki}, 6:65--82, 1960.

\end{thebibliography}

\end{document}